%% file: nips2015.tex
\documentclass{article} 
\usepackage{nips15submit_e,times}
\usepackage{hyperref}
\usepackage{url}

\usepackage{wrapfig}

\usepackage{epsfig}
\usepackage{graphicx}
\usepackage{amsmath}
\usepackage{amssymb}

\usepackage{graphicx} 

\usepackage{overpic}

\usepackage{pgf}

\usepackage{pgfplots}
\newlength\figureheight 
\newlength\figurewidth 

\usepackage{tikz,tkz-euclide,tkz-graph}
\usetikzlibrary{calc}
\usetkzobj{all}

\usepackage{xypic}

\usepackage[noline]{algorithm2e}

\title{MADMM: a generic algorithm for non-smooth optimization on manifolds}

\author{
Artiom Kovnatsky\\
Faculty of Informatics\\
University of Lugano\\ Switzerland\\
\And
Klaus Glashoff\\
Faculty of Informatics\\
University of Lugano\\ Switzerland\\
\And
Michael M. Bronstein\\
Faculty of Informatics\\
University of Lugano\\ Switzerland\\
}

\nipsfinalcopy 

\begin{document}

\maketitle

\begin{abstract}
Numerous problems in machine learning are formulated as optimization with manifold constraints.  
In this paper, we propose the {\em Manifold alternating directions method of multipliers} (MADMM), an extension of the classical 
ADMM scheme for manifold-constrained non-smooth optimization problems and show its application to several challenging problems in dimensionality reduction, data analysis, and manifold learning.  
\end{abstract}

\input{intro.tex}

\input{background.tex}

\input{madmm.tex}

\input{results.tex}

\input{conclusions.tex}

\section*{Acknowledgement} 
The work is supported by the ERC Starting grant No. 307047. 

\small
\bibliographystyle{plain}
\bibliography{nips2015}

\end{document}

%% file: intro.tex
\section{Introduction}

A wide range of problems in machine learning, pattern recognition, computer vision, and signal processing is formulated as optimization problems where the variables are constrained to lie on some Riemannian manifold. 
For example, optimization on the {\em Grassman manifold} comes up in multi-view clustering \cite{dong2014clustering} and matrix completion \cite{keshavan2009gradient}. 
Optimization on the {\em Stiefel manifold}  arises in a plethora of applications ranging from classical ones such as eigenvalue problems, assignment problems, and Procrustes problems \cite{ten1977orthogonal}, to more recent ones such as 1-bit compressed sensing \cite{boufounos20081}. 
Problems involving {\em products of Stiefel manifolds} include shape correspondence \cite{kovnatsky2013coupled}, manifold learning \cite{eynard2015pami}, sensor localization \cite{cucuringu2012sensor}, structural biology \cite{cucuringu2012eigenvector}, and structure from motion recovery \cite{arie2012global}. 
%
%
Optimization on the {\em sphere} is used in principle geodesic analysis \cite{zhang2013probabilistic}, a generalization of the classical PCA to non-Euclidean domains. 
Optimization over the manifold of {\em fixed-rank matrices} arises in 
maxcut problems \cite{journee2010low}, sparse principal component analysis \cite{journee2010low}, regression \cite{meyer2011linear}, matrix completion \cite{boumal2011rtrmc,tan2014riemannian}, and image classification \cite{shalit2010online}. 
{\em Oblique manifolds} 
are encountered in problems such as independent component analysis and joint diagonalization \cite{absil2006joint}, blind source separation \cite{kleinsteuber2012blind}, and 
prediction of stock returns \cite{higham2002computing}. 

%

Though some instances of manifold optimization such as eigenvalues problems have been treated extensively in the distant past, the first general purpose algorithms appeared only in the 1990s 
\cite{smith1994optimization,wen2013feasible}. 
With the emergence of numerous applications during the last decade, especially in the machine learning community,  there has been an increased interest in general-purpose optimization on different manifolds  \cite{absil2009optimization},  leading to several manifold optimization algorithms such as conjugate gradients \cite{edelman1998geometry}, trust regions \cite{absil2007trust}, and Newton \cite{smith1994optimization,alvarez2008unifying}.  
Boumal et al. \cite{boumal2014manopt} released the MATLAB package Manopt, as of today the most complete  generic toolbox for smooth optimization on manifolds, including a variety of manifolds and solvers.

In this paper, we are interested in manifold-constrained minimization of {\em non-smooth} functions. 
%
Recent applications of such problems include several formulations of robust PCA \cite{candes2011robust}, the computation of compressed modes in Euclidean \cite{ozolicnvs2013compressed} and non-Euclidean \cite{neumann2014compressed} spaces, robust Euclidean embedding \cite{cayton2006robust}, synchronization of rotation matrices \cite{wang2013exact}, and functional correspondence 
\cite{huang2014functional,kovnatsky2014functional}.

\paragraph*{Prior work}
Broadly speaking, optimization methods for non-smooth functions break into three classes of approaches. First, 
{\em smoothing} methods replace the non-differentiable objective function with its smooth approximation \cite{chen2012smoothing}. Such methods typically suffer from a tradeoff between accuracy (how far is the smooth approximation from the original objective) and convergence speed (less smooth functions are usually harder to optimize). 
A second class of methods use {\em subgradients} as a generalization of derivatives of non-differentiable functions. In the context of manifold optimization, several subgradient approaches have been proposed \cite{ferreira1998subgradient,ledyaev2007nonsmooth,kleinsteuber2012blind,borckmans2014riemannian}. 
The third class of methods are {\em splitting} approaches. Such methods have been studied mostly for problems involving the minimization of matrix functions with orthogonality constraints. 
%
Lai and Osher proposed the method of splitting orthogonal constraints (SOC) based on the Bregman iteration \cite{lai2014splitting}. Neumann et al. \cite{neumann2014compressed} used a different splitting scheme for the same class of problems.   
%
%

\paragraph*{Contributions}
In this paper, we propose {\em Manifold alternating direction method of multipliers} (MADMM), an extension of the classical 
ADMM scheme \cite{gabay1976dual} for manifold-constrained non-smooth optimization problems. 
The core idea is a splitting into a smooth problem with manifold constraints 
and a non-smooth unconstrained optimization problem. 
Our method has a number of advantages common to ADMM approaches. First, it is very simple to grasp and implement. 
Second, it is generic and not limited to a specific manifold, as opposed e.g. to \cite{lai2014splitting,neumann2014compressed} developed for the Stiefel manifold, or \cite{kleinsteuber2012blind} developed for the oblique manifold. 
Third, it makes very few assumptions about the properties of the objective function.
Fourth,  in some settings, our method lends itself to parallelization on distributed computational architectures \cite{Boyd2010}.  Finally, our method  demonstrates  faster convergence than previous methods in a broad range of applications. 
%
%

%% file: background.tex
\section{Manifold optimization}

The term {\em manifold-} or {\em manifold-constrained optimization} refers to a class of problems of the form 
\begin{eqnarray}
\min_{X \in \mathcal{M}} f(X), 
\label{eq:prob0}
\end{eqnarray}
where $f$ is a smooth real-valued function, $X$ is an $m\times n$ real matrix, and $\mathcal{M}$ is some Riemannian submanifold of $\mathbb{R}^{m\times n}$. 
The manifold is not a vector space and has no global system of coordinates, however, locally at point $X$, the manifold is homeomorphic to a Euclidean space referred to as the {\em tangent space} $T_X \mathcal{M}$.

The main idea of manifold optimization is to treat the objective as a function $f:\mathcal{M} \rightarrow \mathbb{R}$  defined on the manifold,  and perform descent on the manifold itself rather than in the ambient Euclidean space. 
On a manifold, the {\em intrinsic} (Riemannian) gradient $\nabla_{\mathcal{M}}f(X)$ of $f$ at point $X$ is a vector in the tangent space $T_{X}\mathcal{M}$ that can be obtained by projecting the standard (Euclidean) gradient  $\nabla f(X)$ onto $T_{X}\mathcal{M}$ by means of a {\em projection} operator $P_X$ (see an illustration below). 
A step along the intrinsic gradient direction is performed in the tangent plane. In order to obtain the next iterate, the point in the tangent plane is mapped back to the manifold by means of a {\em retraction} operator $R_X$, which is typically an approximation of the {\em exponential map}. 
For many manifolds, the projection $P$ and retraction $R$ operators have a closed form expression.

A conceptual gradient descent-like manifold optimization is presented in Algorithm~\ref{algo:manopt}. For a comprehensive introduction to manifold optimization, the reader is referred to \cite{absil2009optimization}. 



\begin{minipage}{0.6\textwidth}

\RestyleAlgo{boxed}
\IncMargin{0.75em}
\begin{algorithm}[H]

\Repeat{convergence}{
Compute the extrinsic gradient $\nabla f(X^{(k)})$\\
Projection: 
$\nabla_{\mathcal{M}} f(X^{(k)}) = P_{X^{(k)}} (\nabla f(X^{(k)}))$\\
Compute the step size $\alpha^{(k)}$ along the descent direction\\
Retraction: $X^{(k+1)} = R_{X^{(k)}}(- \alpha^{(k)} \nabla_{\mathcal{M}} f(X^{(k)}))$
}
\caption{Conceptual algorithm for smooth optimization on manifold $\mathcal{M}$.}
\label{algo:manopt}
\end{algorithm}

\end{minipage}
\hspace{0.3cm}
\begin{minipage}{0.3\textwidth}
\begin{overpic}
	[width=0.95\linewidth]{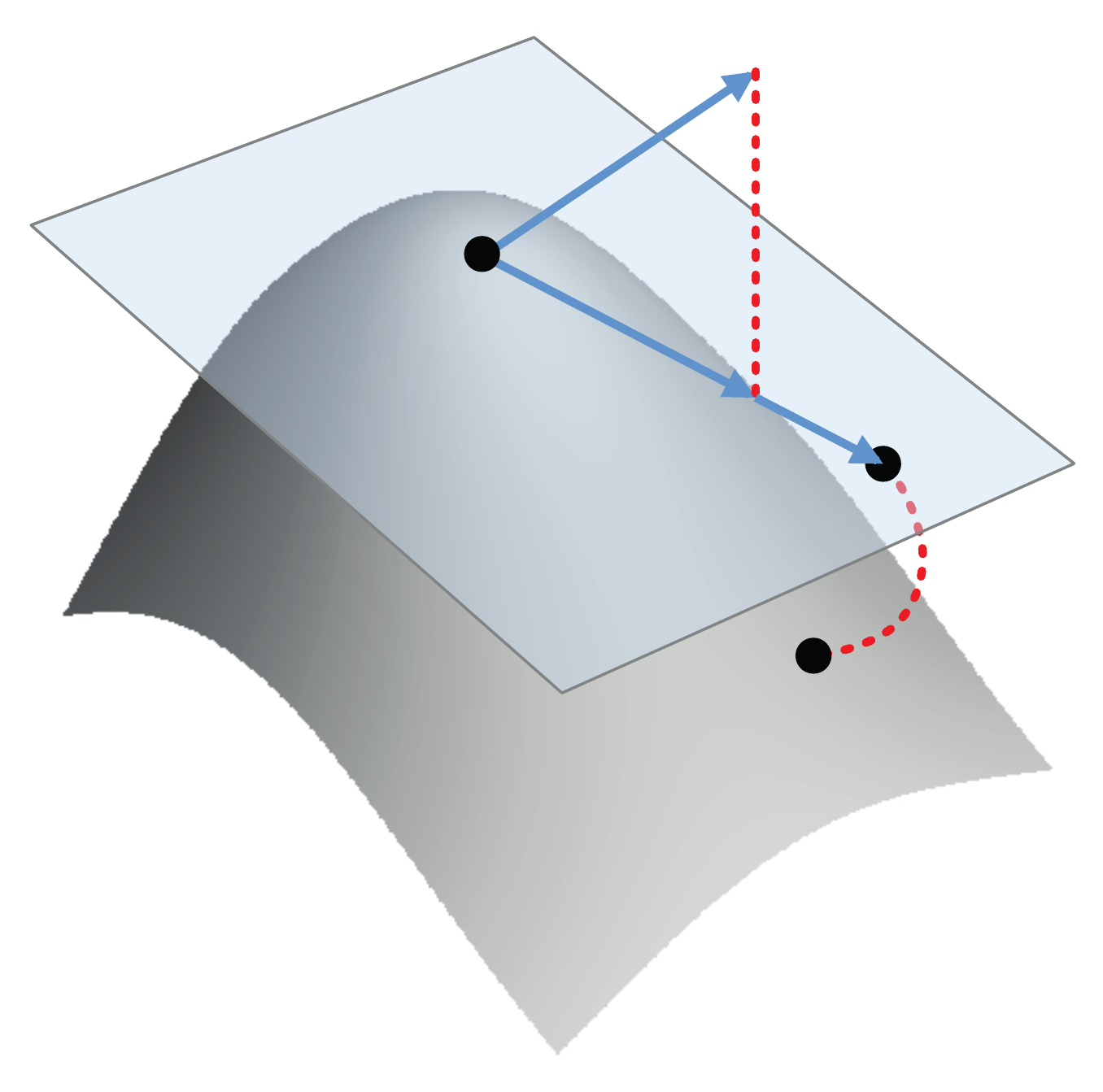}
	\put(39,69){\fontsize{6}{8} $X^{(k)}$}	
	\put(60,95){\fontsize{6}{8} $\nabla f(X^{(k)})$}	
	\put(69,76){\fontsize{6}{8} $P_{X^{(k)}}$}	
	\put(72,60){\fontsize{6}{8} $\alpha^{(k)}\nabla_{\mathcal{M}}f({X^{(k)}})$}	
	\put(83,46){\fontsize{6}{8} $R_{X^{(k)}}$}	
	\put(69,33){\fontsize{6}{8} $X^{(k+1)}$}	
	\put(47,15){\fontsize{6}{8} $\mathcal{M}$}	
\end{overpic}\vspace{0mm}

\end{minipage}

%% file: madmm.tex
\section{Manifold ADMM}

Let us now consider general problems of the form 
\begin{eqnarray}
\min_{X \in \mathcal{M}} f(X) + g(AX), 
\label{eq:prob1}
\end{eqnarray}
where $f$ and $g$ are smooth and non-smooth real-valued functions, respectively, $A$ is a $k\times m$  matrix, and the rest of the notation is as in problem~(\ref{eq:prob0}). 
Examples of $g$ often used in machine learning applications are nuclear-, $L_1$-, or $L_{2,1}$-norms.  
Because of non-smoothness of the objective function, Algorithm~\ref{algo:manopt} cannot be used directly to minimize~(\ref{eq:prob0}).

In this paper, we propose treating this class of problems using the alternating directions method of multipliers (ADMM). 
The key idea is that problem~(\ref{eq:prob1}) can be equivalently formulated as 
\begin{eqnarray}
\min_{X \in \mathcal{M}, Z\in\mathbb{R}^{k\times n}} f(X) + g(Z) \,\,\,\,\, \mathrm{s.t.} \,\,\,\,\, Z = AX
\label{eq:madmm1}
\end{eqnarray}
by introducing an artificial variable $Z$ and a linear constraint.
The method of multipliers \cite{Hestenes1969,Powell1969a}, applied to only the linear constraints in~(\ref{eq:madmm1}), leads to the minimization problem
\begin{eqnarray}
\min_{X \in \mathcal{M}, Z\in\mathbb{R}^{k\times n}} f(X) + g(Z)+\tfrac{\rho }{2}\|AX-Z+U\| _{\mathrm{F}}^2
\label{eq:madmm2}
\end{eqnarray}
where $\rho>0$ and $U\in\mathbb{R}^{k\times n}$ have to be chosen and updated appropriately (see below). 
This formulation now allows splitting the problem into two optimization sub-problems w.r.t. to $X$ and $Z$, which are solved in an alternating manner, followed by an updating of $U$ and, if necessary, of $\rho$. 
Observe that in the first sub-problem w.r.t. $X$ we minimize a {\em smooth} function with manifold constraints, and in the second sub-problem w.r.t. $Z$ we minimize a non-smooth function without manifold constraints. Thus, the problem breaks down into two well-known sub-problems. 
This method, which we call {\em Manifold alternating direction method of multipliers} (MADMM), is summarized in Algorithm~\ref{algo:madmm}.

\RestyleAlgo{boxed}
\IncMargin{0.75em}
\begin{algorithm}[H]
Initialize $k \leftarrow 1$, $Z^{(1)}=X^{(1)}$, $U^{(1)} = 0$. \\
\Repeat{convergence}{
$X$-step: $X^{(k+1)} = \mathop{\mathrm{argmin}}_{X\in \mathcal{M}}  f(X) + \tfrac{\rho }{2}\| AX - Z^{(k)} + U^{(k)}\|_{\mathrm{F}}^2$ \\
$Z$-step: $Z^{(k+1)} = \mathop{\mathrm{argmin}}_{Z}  g(Z) + \tfrac{\rho }{2}\|A X^{(k+1)} - Z + U^{(k)}\|_{\mathrm{F}}^2$\\
Update $U^{(k+1)} = U^{(k)} + AX^{(k+1)} - Z^{(k+1)}$ and $k \leftarrow k+1$\\
}
\caption{Generic MADMM method for non-smooth optimization on manifold $\mathcal{M}$.}
\label{algo:madmm}
\end{algorithm}

Note that MADMM is extremely simple and easy to implement. 
The $X$-step is the setting of Algorithm~\ref{algo:manopt} 
and can be carried out using any standard smooth manifold optimization method. 
Similarly to common implementation of ADMM algorithms, there is no need to solve the $X$-step problem {\em exactly}; instead, only a few iterations of manifold optimization are done.  
Furthermore, for some manifolds and some functions $f$, the $X$-step has a closed-form solution. 
The implementation of the $Z$-step depends on the non-smooth function $g$, and in many cases has a closed-form expression: for example, when $g$ is the $L_1$-norm, the $Z$-step boils down to simple shrinkage, and when $g$ is nuclear norm, the $Z$-step is performed by singular value shrinkage. 
$\rho0$ is the only parameter of the algorithm and its choice is not critical for convergence. In our experiments, we used a rather arbitrary fixed value of $\rho$, though in the ADMM literature it is common to adapt $\rho$ at each iteration, e.g. using the strategy described in \cite{Boyd2010}. 

%% file: results.tex
\section{Results and Applications}

In this section, we show experimental results providing a numerical evaluation of our approach on several challenging applications from the domains of dimensionality reduction, data analysis, pattern recognition, and manifold learning. 
All our experiments were implemented in MATLAB; we used the conjugate gradients and trust regions solvers from the Manopt toolbox \cite{boumal2014manopt} for the $X$-step.  
Time measurements were carried out on a PC with Intel Xeon 2.4 GHz CPU. 

%
%

\subsection{Compressed modes}

\paragraph*{Problem setting}
Our first application is the computation of compressed modes, an approach for constructing localized Fourier-like orthonormal bases recently introduced in \cite{ozolicnvs2013compressed}.  
Let us be given a manifold $\mathcal{S}$ with a Laplacian $\Delta$, where  
in this context, `manifold' can refer to both continuous or discretized manifolds of any dimension, represented as graphs, triangular meshes, etc., and should not be confused with the matrix manifolds we have discussed so far referring to manifold-constrained optimization problems. 
Here, we assume that the manifold is sampled at $n$ points and the Laplacian is represented as an $n\times n$ sparse symmetric matrix. 
In many machine learning applications such as spectral clustering \cite{ng2002spectral}, non-linear dimensionality reduction, and manifold learning \cite{belkin2001laplacian}, one is interested in finding the first $k$ eigenvectors of the Laplacian $\Delta \Phi = \Phi \Lambda$, where $\Phi$ is the $n \times k$ matrix of the first eigenvectors arranged as columns, and $\Lambda$ is the diagonal $k\times k$ matrix of the corresponding eigenvalues.

The first $k$ eigenvectors of the Laplacian can be computed by minimizing the {\em Dirichlet energy} 
\begin{eqnarray}
\min_{\Phi \in \mathbb{R}^{n\times k} } \,\,\, \mathrm{trace}(\Phi^\top \Delta \Phi)  \,\,\,\,\, \mathrm{s.t.} \,\,\,\,\, \Phi^\top \Phi = I 
\label{eq:lap}
\end{eqnarray}
with orthonormality constraints. Laplacian eigenfunctions form an orthonormal basis on the Hilbert space $L^2(\mathcal{S})$ with the standard inner product, and are a generalization of the Fourier basis to non-Euclidean domains. The main disadvantage of such bases is that its elements are globally supported. 
Ozoli{\c{n}}{\v{s}} et al. \cite{ozolicnvs2013compressed} proposed a construction of localized quasi-eigenbases by solving 
\begin{eqnarray}
\min_{ \Phi \in \mathbb{R}^{n\times k}} \,\,\, \mathrm{trace}(\Phi^\top \Delta \Phi) + \mu \| \Phi\|_1 \,\,\,\,\, \mathrm{s.t.} \,\,\,\,\, \Phi^\top \Phi = I,  
\label{eq:cmm}
\end{eqnarray}
where $\mu>0$ is a parameter. 
The $L_1$-norm (inducing sparsity of the resulting basis) together with the Dirichlet energy (imposing smoothness of the basis functions) lead to orthogonal basis functions, referred to as {\em compressed modes} that are localized and approximately diagonalize $\Delta$. 
Lai and Osher \cite{lai2014splitting} and Neumann et al. \cite{neumann2014compressed} proposed two different splitting methods for solving problem~(\ref{eq:cmm}). 
Due to lack of space, the reader is referred to \cite{lai2014splitting,neumann2014compressed} for a detailed description of both methods. 

\paragraph*{Solution}
Here, we realize that problem~(\ref{eq:cmm}) is an instance of manifold optimization on the Stiefel manifold $\mathbb{S}(n,k) = \{ X \in \mathbb{R}^{n\times k} : X^\top X = I \}$ and solve it using MADMM, which assumes in this setting the form of Algorithm~\ref{algo:madmm-cmm}.   
The $X$-step involves optimization of a smooth function on the Stiefel manifold and can be carried out using standard manifold optimization algorithms; we use conjugate gradients and trust regions solvers.  
The $Z$-step requires the minimization of the sum of $L_1$- and $L_2$-norms, a standard problem in signal processing that has an explicit solution by means of thresholding (using the shrinking operator).  
In all our experiments, we used the parameter $\rho = 2$ for MADMM. 
For comparison with the method of \cite{lai2014splitting}, we used the code provided by the authors. The method of \cite{neumann2014compressed} we implemented ourselves.  All the methods were initialized by the same random orthonormal $n\times k$ matrix. 

\RestyleAlgo{boxed}
\begin{algorithm}[H]
\textbf{Input} $n\times n$ Laplacian matrix $\Delta$, parameter $\mu >0$\\
\textbf{Output} $k$ first compressed modes of $\Delta$\\ 
Initialize $k \leftarrow 1$, $\Phi^{(1)} \leftarrow$random orthonormal matrix, $Z^{(1)}=\Phi^{(1)}$, $U^{(1)} = 0$. \\
\Repeat{convergence}{
$\Phi^{(k+1)} = \mathop{\mathrm{argmin}}_{\Phi\in \mathbb{S}(n,k)}  \mathrm{trace}(\Phi^\top \Delta \Phi) + \tfrac{\rho }{2}\| \Phi - Z^{(k)} + U^{(k)}\|_{\mathrm{F}}^2$ \\
%
$Z^{(k+1)} = \mathrm{sign}(\Phi^{(k+1)} + U^{(k)}) \max\{0, |\Phi^{(k+1)} + U^{(k)}| - \tfrac{\mu}{\rho} \}$\\
%
$U^{(k+1)} = U^{(k)} + \Phi^{(k+1)} - Z^{(k+1)} $ and $k \leftarrow k+1$\\
}
\caption{MADMM method for computing compressed modes.}
\label{algo:madmm-cmm}
\end{algorithm}

\begin{figure*}[t!]
\begin{minipage}{1\textwidth}

\centering
\begin{overpic}
	[grid=off,tics=10,trim=0cm 0cm 0cm 0cm,clip,width=1\linewidth]{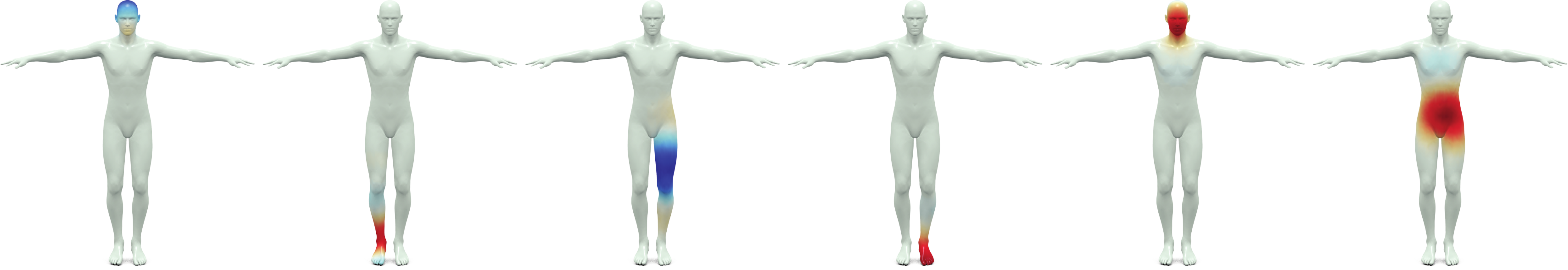}
	%
\end{overpic}\vspace{0mm}
\end{minipage}

\caption{
The first six compressed modes of a human mesh containing $n=8$K points computed using MADMM. Parameter $\mu=10^{-3}$ and three manifold optimization iterations in $X$-step were used in this experiment.   
}
\label{fig:cmms}
\end{figure*}

\begin{figure*}[t!]

\centering
\hspace{-0.75cm}
\begin{minipage}{0.45\textwidth}
\setlength\figurewidth{0.9\linewidth}
\setlength\figureheight{0.6\linewidth} 
\input{cmm_random.tikz}
\end{minipage}
\hspace{0.75cm}
\begin{minipage}{0.45\textwidth}
\setlength\figurewidth{0.9\linewidth}
\setlength\figureheight{0.6\linewidth} 
\input{cmm_solvers.tikz}
\end{minipage}\\
\hspace{-0.75cm}
\begin{minipage}{0.45\textwidth}
\setlength\figurewidth{0.9\linewidth}
\setlength\figureheight{0.6\linewidth} 
\input{cmm_scaling.tikz}
\end{minipage}
\hspace{0.75cm}
\begin{minipage}{0.45\textwidth}
\setlength\figurewidth{0.9\linewidth}
\setlength\figureheight{0.6\linewidth} 
\input{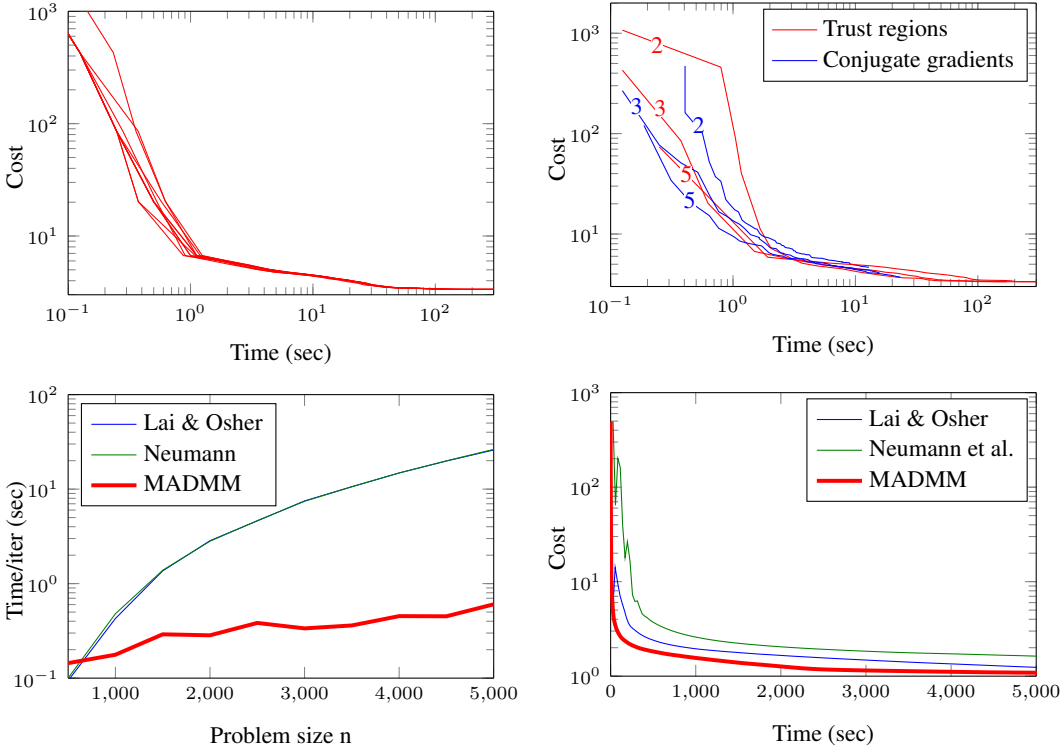}
\end{minipage}

\caption{
Compressed modes problem. 
Top left: convergence of MADMM on a problem of size $n=500, k=10$  
with different random initialization. 
Top right: convergence of MADMM 
using different solvers and number of iterations at $X$-step on the same problem. 
Bottom left: scalability of different methods; shown is time/iteration on a problem of different size (fixed $k=10$ and varying $n$). 
Bottom right: comparison of convergence of different splitting methods and MADMM on a problem of size $n=8$K. 
}
\label{fig:cmm-conv}
\end{figure*}

\paragraph*{Results}
To study the behavior of ADMM, we used a simple 1D problem with a Euclidean Laplacian constructed on a line graph with $n$ vertices.  
Figure~\ref{fig:cmm-conv} (top left) shows the convergence of MADMM with different random initializations. We observe  that the method converges globally independently of the initialization. 
Figure~\ref{fig:cmm-conv} (top right) shows the convergence of MADMM using different solvers and number of iterations in the $X$-step. We did not observe any significant change in the behavior. 
Figure~\ref{fig:cmm-conv} (bottom left) studies the scalability of different algorithms, speaking clearly in favor of MADMM compared to the methods of \cite{lai2014splitting,neumann2014compressed}.
Figure~\ref{fig:cmms} shows compressed modes computed on a triangular mesh of a human sampled at 8K vertices, using the cotangent weights formula \cite{pinkall1993computing} to discretize the Laplacian. 
%
Figure~\ref{fig:cmm-conv} (bottom right) shows the convergence of different methods on this dataset. MADMM shows the best performance among the compared methods.

\subsection{Functional correspondence} 


\paragraph*{Problem setting} 
Our second problem is coupled diagonalization, which is used for finding functional correspondence between manifolds \cite{ovsjanikov2012functional} and multi-view clustering \cite{eynard2015pami}. 
%
%
Let us consider a collection of $L$ manifolds $\{\mathcal{S}_i\}_{i=1}^L$, each discretized at $n_i$ points and equipped with a Laplacian $\Delta_i$ represented as an $n_i \times n_i$ matrix.  
%
%
The {\em functional correspondence} between manifolds $\mathcal{S}_i$ and $\mathcal{S}_j$ is an $n_j\times n_i$ matrix $T_{ij}$ mapping functions from $L^2(\mathcal{S}_i)$ to $L^2(\mathcal{S}_j)$. 
It can be efficiently approximated using the first $k$ Laplacian eigenvectors as $T_{ij} \approx \Phi_j X_{ij} \Phi_i^\top$, where $X_{ij}$ is the $k\times k$ matrix translating Fourier coefficients from basis $\Phi_i$ to basis $\Phi_j$, represented as $n_i\times k$ and $n_j \times k$ matrices, respectively.   
Imposing a further assumption that $T_{ij}$ is volume-preserving, $X_{ij}$ must be an orthonormal matrix \cite{ovsjanikov2012functional}, and thus can be represented as a product of two orthonormal matrices $X_{ij} =  X_i X_j^\top$. 
For each pair of manifolds $\mathcal{S}_i, \mathcal{S}_j$, we assume to be given a set of $q_{ij}$ functions in $L^2(\mathcal{S}_i)$ 
arranged as columns of an $n_i \times q_{ij}$ matrix $F_{ij}$  and the corresponding functions in $L^2(\mathcal{S}_j)$ represented by the $n_j\times q_{ij}$ matrix $G_{ij}$. 
The correspondence between all the manifolds can be established by solving the problem 
\begin{eqnarray}
\min_{X_1, \hdots, X_L} \,\,\, \sum_{i \neq j} \| F_{ij}^\top \Phi_i X_i - G_{ij}^\top \Phi_j X_j  \|_{2,1} + \mu \sum_{i=1}^L \mathrm{off}(X_i^\top \Lambda_i X_i) \,\,\,\,\, \mathrm{s.t.} \,\,\,\,\, X_i^\top X_i = I \,\,\, \forall i=1,\hdots, L, 
\label{eq:funcorr}
\end{eqnarray}
%
%
where $\mathrm{off}(A) = \sum_{i \neq j} a_{ij}^2$, and the use of the $L_{2,1}$-norm 
$\| A \|_{2,1} = \sum_{j} \left( \sum_{i} a_{ij}^2 \right)^{1/2}$
allows to cope with outliers in the correspondence data \cite{wang2013exact,huang2014functional}.  The problem can be interpreted as simultaneous diagonalization of the Laplacians $\Delta_1, \hdots, \Delta_L$ \cite{eynard2015pami}.
As correspondence data $F, G$, one can use point-wise correspondence between some known `seeds', or, in computer graphics applications, some shape descriptors \cite{ovsjanikov2012functional}. 
Geometrically, the matrices $X_i$ can be interpreted as rotations of the respective bases, and the problem tries to achieve a coupling between the bases $\hat{\Phi}_i = \Phi_i X_i$ while making sure that they approximately diagonalize the respective Laplacians. 


\paragraph*{Solution}
Here, we consider problem~(\ref{eq:funcorr}) as optimization on a product of $L$ Stiefel manifolds, $(X_1, \hdots, X_L) \in  \mathbb{S}^L(k,k)$ and solve it using the MADMM method.    
The $X$-step of MADMM in our experiments was performed using four iterations of the manifold conjugate gradients solver. As in the previous problem, the $Z$-step boils down to simple shrinkage. We used $\rho = 1$ and initialized all $X_i = I$. 

\begin{figure*}[t!]

\centering
\hspace{-0.75cm}
\begin{minipage}{0.45\textwidth}
\setlength\figurewidth{0.9\linewidth}
\setlength\figureheight{0.6\linewidth} 
\input{corresp-tosca-kimscurves.tikz}
\end{minipage}
\hspace{0.75cm}
\begin{minipage}{0.45\textwidth}
\setlength\figurewidth{0.9\linewidth}
\setlength\figureheight{0.6\linewidth} 
\input{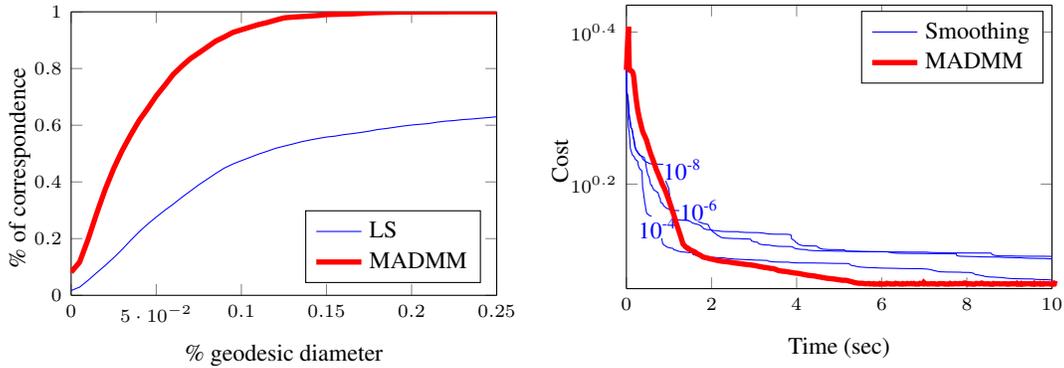}
\end{minipage}

\caption{
Functional correspondence problem. 
Left: evaluation of the functional correspondence obtained using MADMM and least squares. 
Shown in the percentage of correspondences falling within a geodesic ball of increasing radius w.r.t. the groundtruth correspondence. 
Right: convergence of MADMM and smoothing method for various values of the smoothing parameter. 
}
\label{fig:corr-curves}
\end{figure*}

\begin{figure*}[t!]
\begin{minipage}{1\textwidth}

\centering
\begin{overpic}
	[grid=off,tics=10,trim=0cm 0cm 0cm 0cm,clip,width=1\linewidth]{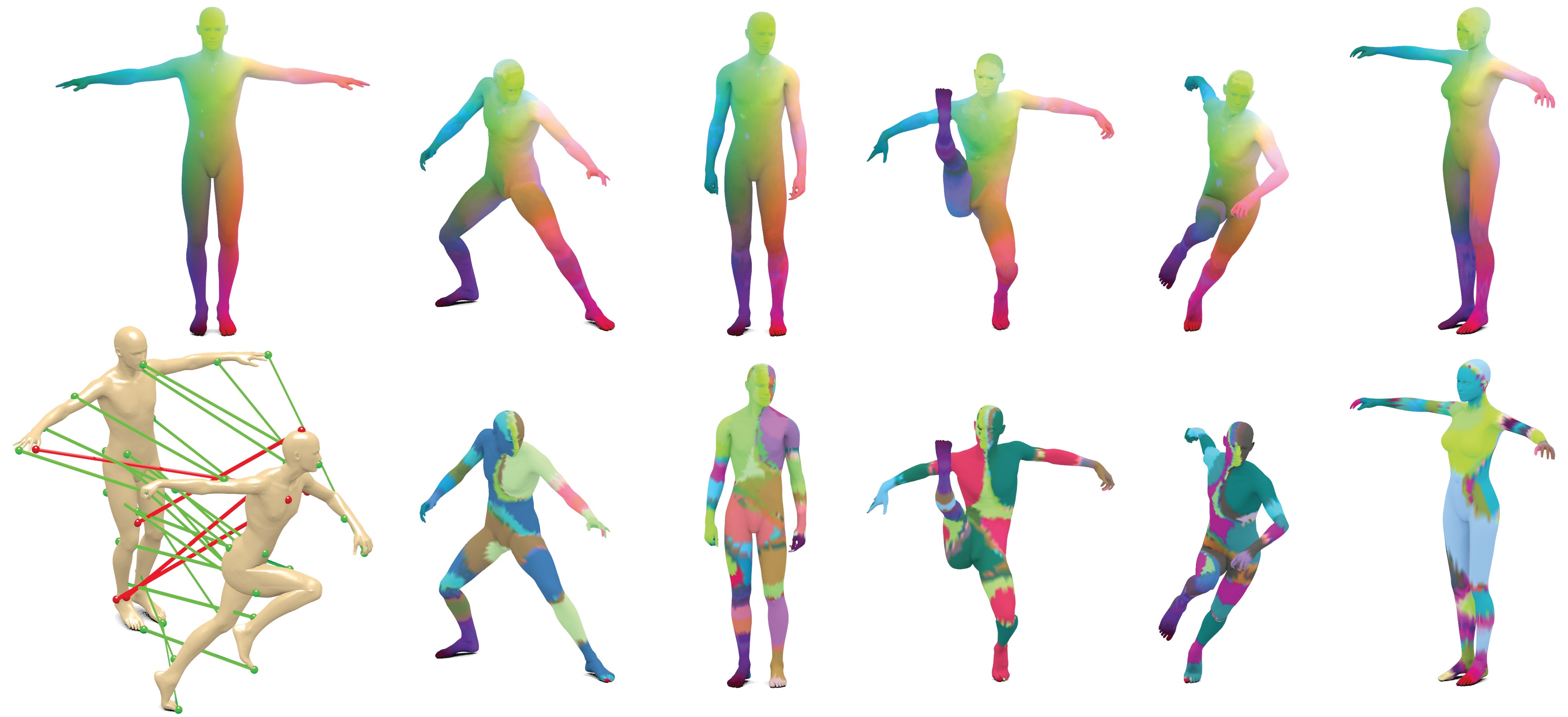}
	%
\end{overpic}\vspace{0mm}
\end{minipage}

\caption{
Examples of correspondences obtained with MADMM (top) and least-squares solution (bottom). Similar colors encode corresponding points. Bottom left: examples of correspondence between a pair of shapes (outliers are shown in red). 
}
\label{fig:corr-shapes}
\end{figure*}

\paragraph*{Results} 
We computed functional correspondences between $L=6$ human 3D shapes from the TOSCA dataset \cite{bronstein2008numerical} 
using  $k=25$ basis functions and $q=25$ seeds as correspondence data, contaminated by 16\% outliers. 
Figure~\ref{fig:corr-curves} (left) analyzes the resulting correspondence quality using the Princeton protocol \cite{kim2011blended}, plotting the percentage of correspondences falling within a geodesic ball of increasing radius w.r.t. the groundtruth correspondence. 
For comparison, we show the results of a least-squares solution used in \cite{ovsjanikov2012functional} (see Figure~\ref{fig:corr-shapes}).  
%
%
Figure~\ref{fig:corr-curves} (right) shows the convergence of MADMM in a correspondence problem with $L=2$ shapes. For comparison, we show the convergence of a smoothed version of the $L_{2,1}$-norm $\| A \|_{2,1} \approx \sum_{j} \left( \sum_{i} a_{ij}^2 + \epsilon \right)^{1/2}$ in~(\ref{eq:funcorr}) for various values of the smoothing parameter $\epsilon$. 

\subsection{Robust Euclidean embedding}

\paragraph*{Problem setting}
Our third problem is an $L_1$ formulation of the multidimensional scaling (MDS) problem treated in \cite{cayton2006robust} under the name {\em robust Euclidean embedding} (REE). 
Let us be given an $n\times n$ matrix $D$ of squared distances. The goal is to find a $k$-dimensional configuration of points $X \in \mathbb{R}^{n\times k}$ such that the Euclidean distances between them are as close as possible to the given ones. 
The classical MDS approach employs the duality between Euclidean distance matrices and Gram matrices: a squared Euclidean distance matrix $D$ can be converted into a {\em similarity matrix} by means of double-centering $B = -\tfrac{1}{2}HDH$, where $H = I - \tfrac{1}{n}11^\top$. 
Conversely, the squared distance matrix is obtained from $B$ by $(\mathrm{dist}(B))_{ij} = b_{ii} + b_{jj} - 2 b_{ij}$. 
The similarity matrix corresponding to a Euclidean distance matrix is positive semi-definite and can be represented as a Gram matrix $B = XX^\top$, where $X$ is the desired embedding.  In the case when $D$ is not Euclidean, $B$ acts as a low-rank approximation of the similarity matrix (now not necessarily positive semi-definite) associated with $D$, leading to the problem  
\begin{eqnarray}
\min_{X \in \mathbb{R}^{m\times k}} \| HDH - XX^\top \|_\mathrm{F}^2
\end{eqnarray}
known as {\em classical MDS} or {\em classical scaling}, which has a closed form solution by means of eigendecomposition of $HDH$.

The main disadvantage of classical MDS is the fact that noise in a single entry of the distance matrix $D$ is spread over entire column/row by the double centering transformation. To cope with this problem, Cayton and Dasgupta \cite{cayton2006robust} proposed an $L_1$ version of the problem, 
\begin{eqnarray}
\min_{B \in \mathbb{R}^{n\times n} } \| D - \mathrm{dist}(B) \|_1 \,\,\,\,\, \mathrm{s.t.} \,\,\,\,\,
B \succeq 0, \, \mathrm{rank}(B) \leq k,  
\label{eq:mds1}
\end{eqnarray}
where the use of the $L_1$-norm efficiently rejects outliers. 
The authors proposed two solutions for problem~(\ref{eq:mds1}): a semi-definite programming (SDP) formulation and a subgradient descent algorithm 
%
(the reader is referred to \cite{cayton2006robust} for a detailed description of both methods).

\paragraph*{Solution}
Here, we consider~(\ref{eq:mds1}) as a non-smooth optimization of the form~(\ref{eq:prob1}) on the manifold of fixed-rank positive semi-definite matrices and solve it using MADMM. Note that in this case, we have only the non-smooth function $g$ and $f \equiv 0$. 
The $X$-step of the MADMM algorithm is manifold optimization of a quadratic function, carried out using two iterations of manifold conjugate gradients solver. 
The $Z$-step is performed by shrinkage.   
In our experiments, all the compared methods were initialized with the classical MDS solution and the value $\rho = 10$ was used for MADMM. SDP approach was implemented using MATLAB CVX toolbox \cite{grant2008cvx}. 

%
%
%

\begin{figure*}[t!]
\centering
\setlength\figureheight{0.5\linewidth} 
\setlength\figurewidth{0.6\linewidth}
\input{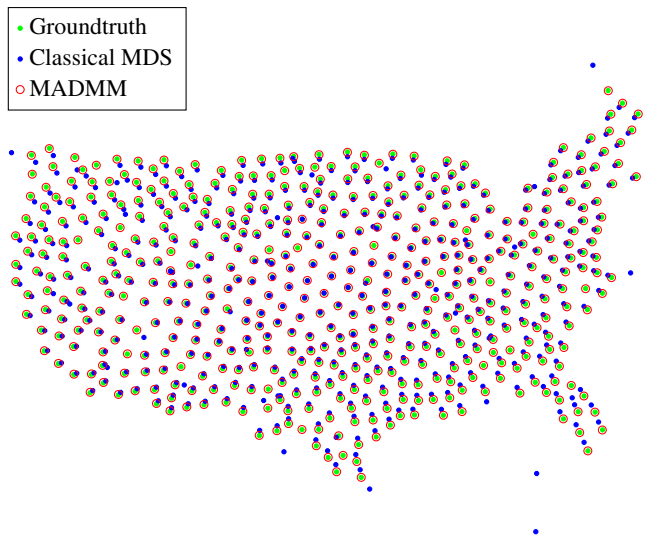}
\caption{
Embedding of the noisy distances between 500 US cities in the plane using classical MDS (blue) and REE solved using MADMM (red). The distance matrix was contaminated by sparse noise by doubling the distance between some cities. Note how classical MDS is sensitive to  outliers. 
}
\label{fig:mds-cities}
\end{figure*}

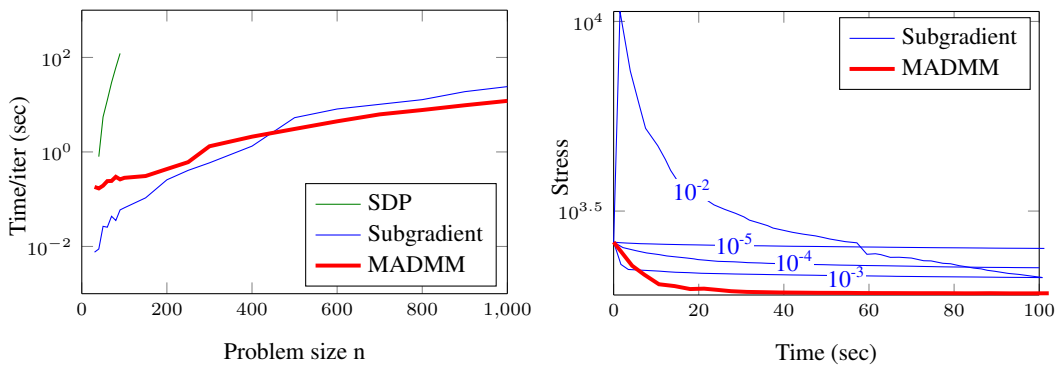
\begin{figure*}

\centering
\hspace{-0.75cm}
\begin{minipage}{0.45\textwidth}
\setlength\figurewidth{0.9\linewidth}
\setlength\figureheight{0.6\linewidth} 
\input{mds_scale.tikz}
\end{minipage}
\hspace{0.75cm}
\begin{minipage}{0.45\textwidth}
\setlength\figurewidth{0.9\linewidth}
\setlength\figureheight{0.6\linewidth} 
\input{mds_convergence.tikz}
\end{minipage}

\caption{
REE problem. 
Left: scalability of different algorithms; shown is single iteration complexity as functions of the problem size $n$ using random distance data. SDP did not scale beyond $n=100$. 
Right: example of convergence of MADMM and subgradient algorithm of \cite{cayton2006robust} on the US cities problem of size $n=500$. The subgradient algorithm is very sensitive to the choice of the initial step size $c$; choosing too large $c$ breaks the convergence.   
}
\label{fig:mds-conv}
\end{figure*}

\paragraph*{Results}
Figure~\ref{fig:mds-cities} shows an example of 2D Euclidean embedding of the distances between 500 US cities, contaminated by sparse noise (the distance between two cities was doubled, as in a similar experiment in \cite{cayton2006robust}). The robust embedding is insensitive to such outliers, while the classical MDS result is completely ruined. 
Figure~\ref{fig:mds-conv} (right) shows an example of convergence of the proposed MADMM method and the subgradient descent of \cite{cayton2006robust} on the same dataset. 
We observed that our algorithm outperforms the subgradient method in terms of convergence speed.  
Furthermore, the subgradient method appears to be very sensitive to the initial step size $c$; choosing too small a step leads to slower convergence, and if the step is too large the algorithm may fail to converge. 
Figure~\ref{fig:mds-conv} (left) studies the scalability of the subgradient-, SDP-, and MADMM-based  solutions for the REE problem, plotting the complexity of a single iteration as function of the problem size on random data. The typical number of iterations was of the order of 20 for SDP, 50 for MADMM, and 500 for the subgradient method. We see that MADMM scaled better than other approaches, and that SDP is not applicable to large problems.


%% file: cmm_random.tikz
%
%
\begin{tikzpicture}
\tikzstyle{every node}=[font=\footnotesize]

\begin{axis}[%
width=\figurewidth,
height=\figureheight,
at={(0\figurewidth,0\figureheight)},
scale only axis,
separate axis lines,
every outer x axis line/.append style={black},
every x tick label/.append style={font=\color{black}},
every x tick label/.append style={font=\scriptsize},
xmode=log,
xmin=0.1,
xmax=300,
xminorticks=true,
xlabel={Time (sec)},
every outer y axis line/.append style={black},
every y tick label/.append style={font=\color{black}},
every y tick label/.append style={font=\scriptsize},
ymode=log,
ymin=3,
ymax=1000,
yminorticks=true,
y label style={at={(axis description cs:0.05,0.45)},anchor=north},
ylabel={Cost}
]
\addplot [color=red,solid,forget plot]
  table[row sep=crcr]{%
0.136083697331127	1099.7\\
0.23400150006637	428.341302934067\\
0.358802300062962	85.1774040645901\\
0.624003999982961	20.2496617476092\\
1.24800799996592	6.70323948650063\\
2.74561760004144	5.55515768956423\\
4.11842640000395	5.06002491404293\\
5.86563760007266	4.77667373051705\\
8.23685280000791	4.57707006764094\\
11.1072712000459	4.42507464956054\\
13.4784863999812	4.29311535638748\\
16.099303200026	4.15732466205045\\
18.9697216000641	4.06914434470653\\
21.3409367999993	3.96119206446315\\
23.712152000051	3.87585558717682\\
26.2081679999828	3.79275988145233\\
28.4545824000379	3.6866923278836\\
30.57619599998	3.63038038261034\\
33.446614400018	3.59193022974517\\
36.5666344000492	3.55112138756421\\
39.1874511999777	3.50783970120968\\
42.9938756000483	3.48705832576353\\
46.6130988000659	3.47150663045471\\
50.0451208000304	3.45924198930155\\
54.1635472000344	3.44920263184368\\
57.6579696000554	3.44078731625236\\
60.9027903999668	3.43421246261438\\
64.022810399998	3.43070871410766\\
67.2676312000258	3.42738381375213\\
70.7620536000468	3.42115586500824\\
74.6308784000576	3.41520282136783\\
78.4061026000418	3.41137291050434\\
82.3997282000491	3.40750202803593\\
86.1437522000633	3.40377936476731\\
89.2637721999781	3.40011304407417\\
92.9453958000522	3.39617362238874\\
96.3150173999602	3.39203246280603\\
99.8094397999812	3.38768107881226\\
103.678264599992	3.38303757735749\\
107.188287100056	3.37824903955962\\
110.557908699964	3.37348068468192\\
114.052331099985	3.36890816472437\\
118.342358599999	3.36467269715585\\
121.46237860003	3.36115264062007\\
124.582398600061	3.35828388405347\\
128.076820999966	3.35592761656574\\
132.070446599973	3.35423740657108\\
136.313673799974	3.35285813724321\\
140.556900999974	3.35185901633073\\
144.550526599982	3.35098088372188\\
149.168156199972	3.35034071368936\\
154.035387400072	3.3498944019755\\
158.465815800009	3.34949957309996\\
162.584242200013	3.34911457627257\\
167.389073000057	3.34887004551873\\
172.505905800033	3.34868390103383\\
178.02834119997	3.34849500493736\\
182.770771600073	3.34831739717445\\
187.762803600053	3.34818070539874\\
190.258819599985	3.34807947588771\\
192.630034800037	3.34801538378178\\
195.126050799969	3.34797005648034\\
197.996469200007	3.34789640321076\\
200.617286000052	3.34782319311267\\
203.737305999966	3.34779772290366\\
206.576524200034	3.3477782040913\\
208.94773939997	3.3477373001841\\
211.443755400018	3.34769779913119\\
213.939771400066	3.34766145911694\\
216.560588199995	3.34764663469974\\
219.431006600033	3.34762281998074\\
221.802221799968	3.34759860410849\\
224.672640200006	3.34759098637478\\
227.66785940004	3.34757245283431\\
230.039074599976	3.34755440688097\\
232.784692200017	3.34752638422135\\
235.155907400069	3.34749690343533\\
237.527122600004	3.34747761847885\\
239.898337800056	3.34746738306767\\
242.519154599984	3.34746378661995\\
245.763975400012	3.34746307258206\\
248.384792200057	3.34746145510143\\
250.880808199989	3.34746027348737\\
253.501625000034	3.34745919893054\\
256.247242600075	3.34744635263903\\
258.868059400003	3.34742750661081\\
261.364075400052	3.34742038546519\\
263.98489219998	3.34741633333558\\
266.605709000025	3.34741749672954\\
269.741329099983	3.34741967663356\\
272.861349100014	3.3474164343982\\
275.482165900059	3.34741854199888\\
278.227783499984	3.34742339788904\\
281.223002700019	3.34741899461979\\
284.218221900053	3.34740859388272\\
287.088640299975	3.34740125878164\\
290.208660300006	3.34738625633781\\
293.079078700044	3.34738357686429\\
296.261499100016	3.34737478632206\\
300.00552310003	3.34736818244682\\
};
\addplot [color=red,solid,forget plot]
  table[row sep=crcr]{%
0.0725779718873673	1099.7\\
0.124800799996592	428.341302934067\\
0.249601599993184	85.1774040645901\\
0.374402399989776	20.2496617476092\\
1.09200699999928	6.70323948650063\\
2.83921819995157	5.55515768956423\\
4.21202699991409	5.06002491404293\\
5.95923819998279	4.77667373051705\\
8.33045339991804	4.57707006764094\\
10.8264693999663	4.42507464956054\\
13.072883799905	4.29311535638748\\
15.3660984999733	4.15732466205045\\
17.8621144999051	4.06914434470653\\
20.3581304999534	3.96119206446315\\
22.9789472999983	3.87585558717682\\
25.9273661999032	3.79275988145233\\
28.548182999948	3.6866923278836\\
31.0441989999963	3.63038038261034\\
33.7898165999213	3.59193022974517\\
36.9098365999525	3.55112138756421\\
40.2794581999769	3.50783970120968\\
43.3994782000082	3.48705832576353\\
46.8939005999127	3.47150663045471\\
50.5131237999303	3.45924198930155\\
54.6315501999343	3.44920263184368\\
58.2507733999519	3.44078731625236\\
61.7451957999729	3.43421246261438\\
64.9900166000007	3.43070871410766\\
68.6092398000183	3.42738381375213\\
72.353263799916	3.42115586500824\\
75.847686199937	3.41520282136783\\
79.342108599958	3.41137291050434\\
82.836530999979	3.40750202803593\\
86.4557541999966	3.40377936476731\\
89.5757741999114	3.40011304407417\\
93.194997399929	3.39617362238874\\
95.6910133999772	3.39203246280603\\
99.1854357999982	3.38768107881226\\
102.929459800012	3.38303757735749\\
106.548682999914	3.37824903955962\\
109.918304599938	3.37348068468192\\
112.913523799973	3.36890816472437\\
116.407946199994	3.36467269715585\\
120.151970200008	3.36115264062007\\
123.396790999919	3.35828388405347\\
127.140814999933	3.35592761656574\\
131.134440599941	3.35423740657108\\
135.377667799941	3.35285813724321\\
139.620894999942	3.35185901633073\\
143.115317399963	3.35098088372188\\
147.920148200006	3.35034071368936\\
153.270982499933	3.3498944019755\\
159.557822800009	3.34949957309996\\
164.300253199996	3.34911457627257\\
170.337491899962	3.34887004551873\\
174.955121499952	3.34868390103383\\
180.196755099925	3.34849500493736\\
185.313587900018	3.34831739717445\\
190.867223499925	3.34818070539874\\
193.86244269996	3.34807947588771\\
195.859255499905	3.34801538378178\\
198.10566989996	3.34797005648034\\
200.851287500001	3.34789640321076\\
203.721705899923	3.34782319311267\\
206.342522699968	3.34779772290366\\
209.088140300009	3.3477782040913\\
211.708957099938	3.3477373001841\\
214.329773899983	3.34769779913119\\
217.090991599951	3.34766145911694\\
219.836609199992	3.34764663469974\\
222.831828399911	3.34762281998074\\
225.827047599945	3.34759860410849\\
228.44786439999	3.34759098637478\\
230.819079599925	3.34757245283431\\
233.315095599974	3.34755440688097\\
235.686310799909	3.34752638422135\\
238.05752599996	3.34749690343533\\
241.052745199995	3.34747761847885\\
244.047964399913	3.34746738306767\\
247.292785199941	3.34746378661995\\
250.537605999969	3.34746307258206\\
253.907227599993	3.34746145510143\\
257.027247599908	3.34746027348737\\
260.022466799943	3.34745919893054\\
262.892885199981	3.34744635263903\\
266.168906199979	3.34742750661081\\
269.17972549994	3.34742038546519\\
272.674147899961	3.34741633333558\\
275.918968699989	3.34741749672954\\
278.539785499917	3.34741967663356\\
281.160602299962	3.3474164343982\\
284.0310207	3.34741854199888\\
287.151040699915	3.34742339788904\\
289.896658299956	3.34741899461979\\
292.891877499991	3.34740859388272\\
296.011897499906	3.34740125878164\\
299.131917499937	3.34738625633781\\
302.127136699972	3.34738357686429\\
};
\addplot [color=red,solid,forget plot]
  table[row sep=crcr]{%
0.0725779718873673	1099.7\\
0.124800799996592	428.341302934067\\
0.249601599993184	85.1774040645901\\
0.499203199986368	20.2496617476092\\
0.998406399972737	6.70323948650063\\
2.62081680004485	5.55515768956423\\
3.99362560000736	5.06002491404293\\
5.8344373999862	4.77667373051705\\
8.58005500002764	4.57707006764094\\
11.2008717999561	4.42507464956054\\
13.197684600018	4.29311535638748\\
15.943302199943	4.15732466205045\\
18.813720599981	4.06914434470653\\
21.3097366000293	3.96119206446315\\
24.4297565999441	3.87585558717682\\
27.3001749999821	3.79275988145233\\
29.5465894000372	3.6866923278836\\
31.7930037999759	3.63038038261034\\
34.5386214000173	3.59193022974517\\
38.0330438000383	3.55112138756421\\
41.6522669999395	3.50783970120968\\
44.8970877999673	3.48705832576353\\
47.5179046000121	3.47150663045471\\
50.2635221999371	3.45924198930155\\
53.1339405999752	3.44920263184368\\
56.6283629999962	3.44078731625236\\
59.873183800024	3.43421246261438\\
63.367606200045	3.43070871410766\\
66.8620285999496	3.42738381375213\\
69.9040480999975	3.42115586500824\\
73.0708684000419	3.41520282136783\\
76.7836921999697	3.41137291050434\\
80.5277161999838	3.40750202803593\\
83.8973378000082	3.40377936476731\\
87.2669594000326	3.40011304407417\\
90.7613817999372	3.39617362238874\\
93.8814017999684	3.39203246280603\\
97.0014217999997	3.38768107881226\\
100.495844200021	3.38303757735749\\
104.879472299945	3.37824903955962\\
108.748297099955	3.37348068468192\\
112.367520299973	3.36890816472437\\
115.487540300004	3.36467269715585\\
118.732361100032	3.36115264062007\\
122.725986700039	3.35828388405347\\
126.345209899941	3.35592761656574\\
130.463636299944	3.35423740657108\\
134.457261899952	3.35285813724321\\
138.466487600002	3.35185901633073\\
142.46011320001	3.35098088372188\\
147.202543599997	3.35034071368936\\
152.319376399973	3.3498944019755\\
157.186607599957	3.34949957309996\\
161.554635599954	3.34911457627257\\
166.421866799938	3.34887004551873\\
171.039496400044	3.34868390103383\\
174.783520399942	3.34849500493736\\
179.525950800045	3.34831739717445\\
184.471182500012	3.34818070539874\\
186.343194499961	3.34807947588771\\
188.46480810002	3.34801538378178\\
190.836023299955	3.34797005648034\\
193.503640400013	3.34789640321076\\
196.623660400044	3.34782319311267\\
199.400478199939	3.34779772290366\\
201.662492700038	3.3477782040913\\
204.439310500049	3.3477373001841\\
206.560924099991	3.34769779913119\\
208.807338500046	3.34766145911694\\
211.927358499961	3.34764663469974\\
214.548175300006	3.34762281998074\\
217.54339450004	3.34759860410849\\
220.289012099965	3.34759098637478\\
223.159430500004	3.34757245283431\\
225.905048100045	3.34755440688097\\
228.65066569997	3.34752638422135\\
231.271482500015	3.34749690343533\\
234.266701700049	3.34747761847885\\
237.386721699964	3.34746738306767\\
240.381940899999	3.34746378661995\\
243.377160100034	3.34746307258206\\
246.497180099948	3.34746145510143\\
249.492399299983	3.34746027348737\\
252.908821200021	3.34745919893054\\
255.404837199952	3.34744635263903\\
258.025653999997	3.34742750661081\\
261.270474800025	3.34742038546519\\
264.64009640005	3.34741633333558\\
267.760116399964	3.34741749672954\\
270.880136399996	3.34741967663356\\
274.000156400027	3.3474164343982\\
277.790980699938	3.34741854199888\\
281.035801499966	3.34742339788904\\
283.65661830001	3.34741899461979\\
285.903032699949	3.34740859388272\\
288.898251899984	3.34740125878164\\
291.519068700029	3.34738625633781\\
294.264686299954	3.34738357686429\\
297.259905499988	3.34737478632206\\
300.27072479995	3.34736818244682\\
};
\addplot [color=red,solid,forget plot]
  table[row sep=crcr]{%
0.0725779718873673	1099.7\\
0.124800799996592	428.341302934067\\
0.249601599993184	85.1774040645901\\
0.499203199986368	20.2496617476092\\
0.998406399972737	6.70323948650063\\
2.74561760004144	5.55515768956423\\
4.24322720000055	5.06002491404293\\
6.11523920006584	4.77667373051705\\
8.61125519999769	4.57707006764094\\
10.6080680000596	4.42507464956054\\
12.7764819000149	4.29311535638748\\
15.1476971000666	4.15732466205045\\
18.1429162999848	4.06914434470653\\
20.8885339000262	3.96119206446315\\
23.3845498999581	3.87585558717682\\
25.7557651000097	3.79275988145233\\
28.1269803000614	3.6866923278836\\
31.2938005999895	3.63038038261034\\
33.5402150000446	3.59193022974517\\
36.0362309999764	3.55112138756421\\
39.7022545000073	3.50783970120968\\
42.6974737000419	3.48705832576353\\
45.4430912999669	3.47150663045471\\
49.0623144999845	3.45924198930155\\
52.6815377000021	3.44920263184368\\
56.5503625000129	3.44078731625236\\
60.0447849000338	3.43421246261438\\
63.2896057000617	3.43070871410766\\
66.2848248999799	3.42738381375213\\
69.6544465000043	3.42115586500824\\
73.0240681000287	3.41520282136783\\
76.2844889999833	3.41137291050434\\
80.3405150000472	3.40750202803593\\
84.5837422000477	3.40377936476731\\
89.2481721000513	3.40011304407417\\
92.9921961000655	3.39617362238874\\
96.3618176999735	3.39203246280603\\
100.105841699988	3.38768107881226\\
103.974666499998	3.38303757735749\\
107.593889700016	3.37824903955962\\
110.713909700047	3.37348068468192\\
114.255132399965	3.36890816472437\\
117.375152399996	3.36467269715585\\
120.245570800034	3.36115264062007\\
123.864794000052	3.35828388405347\\
127.608818000066	3.35592761656574\\
130.604037199984	3.35423740657108\\
134.472861999995	3.35285813724321\\
138.466487600002	3.35185901633073\\
142.725314900046	3.35098088372188\\
147.21814370004	3.35034071368936\\
151.710972500034	3.3498944019755\\
156.32860210014	3.34949957309996\\
161.195833300124	3.34911457627257\\
165.688662100118	3.34887004551873\\
169.807088500122	3.34868390103383\\
174.674319700105	3.34849500493736\\
179.666351700085	3.34831739717445\\
183.784778100089	3.34818070539874\\
186.031192500144	3.34807947588771\\
188.776810100069	3.34801538378178\\
191.647228500107	3.34797005648034\\
194.268045300036	3.34789640321076\\
196.764061300084	3.34782319311267\\
199.384878100129	3.34779772290366\\
201.631292500068	3.3477782040913\\
204.501710900106	3.3477373001841\\
207.372129300144	3.34769779913119\\
210.258147800108	3.34766145911694\\
213.003765400033	3.34764663469974\\
215.874183800071	3.34762281998074\\
218.495000600116	3.34759860410849\\
221.365419000038	3.34759098637478\\
224.485439000069	3.34757245283431\\
227.6054590001	3.34755440688097\\
230.101475000149	3.34752638422135\\
232.472690200084	3.34749690343533\\
235.467909400119	3.34747761847885\\
238.712730200146	3.34746738306767\\
241.957551000058	3.34746378661995\\
245.061970900046	3.34746307258206\\
247.557986900094	3.34746145510143\\
250.475205600145	3.34746027348737\\
253.59522560006	3.34745919893054\\
256.465644000098	3.34744635263903\\
260.256468300126	3.34742750661081\\
263.37648830004	3.34742038546519\\
265.997305100085	3.34741633333558\\
269.05492470006	3.34741749672954\\
272.798948700074	3.34741967663356\\
276.168570300099	3.3474164343982\\
279.46019140014	3.34741854199888\\
282.580211400054	3.34742339788904\\
285.450629800092	3.34741899461979\\
289.61585650011	3.34740859388272\\
292.361474100035	3.34740125878164\\
295.107091700076	3.34738625633781\\
297.852709300118	3.34738357686429\\
300.473526100046	3.34737478632206\\
};
\addplot [color=red,solid,forget plot]
  table[row sep=crcr]{%
0.0725779718873673	1099.7\\
0.124800799996592	428.341302934067\\
0.374402399989776	85.1774040645901\\
0.624003999982961	20.2496617476092\\
1.12320719996933	6.70323948650063\\
2.37121519993525	5.55515768956423\\
3.61922319990117	5.06002491404293\\
5.11683279997669	4.77667373051705\\
7.36324719991535	4.57707006764094\\
10.1088647999568	4.42507464956054\\
12.7296816000016	4.29311535638748\\
15.6000999999233	4.15732466205045\\
18.3457175999647	4.06914434470653\\
20.9665344000096	3.96119206446315\\
24.632557899924	3.87585558717682\\
27.736977799912	3.79275988145233\\
30.7789972999599	3.6866923278836\\
33.8990172999911	3.63038038261034\\
36.5042339999927	3.59193022974517\\
39.7490547999041	3.55112138756421\\
42.9938755999319	3.50783970120968\\
46.3634971999563	3.48705832576353\\
49.6083179999841	3.47150663045471\\
52.9779396000085	3.45924198930155\\
56.8935647000326	3.44920263184368\\
60.5127879000502	3.44078731625236\\
63.8824094999582	3.43421246261438\\
67.3768318999792	3.43070871410766\\
70.9960550999967	3.42738381375213\\
74.7400791000109	3.42115586500824\\
78.1097007000353	3.41520282136783\\
83.0393322999589	3.41137291050434\\
86.9081570999697	3.40750202803593\\
90.2777786999941	3.40377936476731\\
93.7722011000151	3.40011304407417\\
97.2822235999629	3.39617362238874\\
100.527044399991	3.39203246280603\\
103.522263600025	3.38768107881226\\
106.64228359994	3.38303757735749\\
110.261506799958	3.37824903955962\\
113.490727499942	3.37348068468192\\
116.73554829997	3.36890816472437\\
119.605966700008	3.36467269715585\\
123.100389100029	3.36115264062007\\
126.844413100043	3.35828388405347\\
130.838038700051	3.35592761656574\\
134.332461099955	3.35423740657108\\
138.575688299956	3.35285813724321\\
142.943716299953	3.35185901633073\\
146.438138699974	3.35098088372188\\
150.306963499985	3.35034071368936\\
155.595397399971	3.3498944019755\\
159.713823799975	3.34949957309996\\
164.893057000008	3.34911457627257\\
170.009889799985	3.34887004551873\\
175.126722599962	3.34868390103383\\
181.538363700034	3.34849500493736\\
186.389994799974	3.34831739717445\\
191.756429199944	3.34818070539874\\
194.377245999989	3.34807947588771\\
196.873262000037	3.34801538378178\\
199.244477199973	3.34797005648034\\
201.865294000017	3.34789640321076\\
204.361309999949	3.34782319311267\\
206.607724400004	3.34779772290366\\
208.729337999946	3.3477782040913\\
211.225353999995	3.3477373001841\\
213.970971600036	3.34769779913119\\
216.841389999958	3.34766145911694\\
219.212605200009	3.34764663469974\\
222.083023600047	3.34762281998074\\
225.078242799966	3.34759860410849\\
227.823860400007	3.34759098637478\\
230.819079600042	3.34757245283431\\
233.939099599957	3.34755440688097\\
236.934318799991	3.34752638422135\\
240.054338800022	3.34749690343533\\
243.174358799937	3.34747761847885\\
246.294378799968	3.34746738306767\\
249.4143988	3.34746378661995\\
252.534418800031	3.34746307258206\\
255.404837199952	3.34746145510143\\
258.524857199984	3.34746027348737\\
261.894478800008	3.34745919893054\\
264.889698000043	3.34744635263903\\
268.680522299954	3.34742750661081\\
272.424546299968	3.34742038546519\\
275.794167899992	3.34741633333558\\
278.539785500034	3.34741749672954\\
281.535004699952	3.34741967663356\\
284.43662329996	3.3474164343982\\
287.556643299991	3.34741854199888\\
290.427061700029	3.34742339788904\\
293.297480099951	3.34741899461979\\
296.292699299986	3.34740859388272\\
299.28791850002	3.34740125878164\\
302.033536099945	3.34738625633781\\
};
\addplot [color=red,solid,forget plot]
  table[row sep=crcr]{%
0.0725779718873673	1099.7\\
0.124800799996592	428.341302934067\\
0.280801799963228	85.1774040645901\\
0.530403399956413	20.2496617476092\\
0.904805799946189	6.70323948650063\\
2.15281379991211	5.55515768956423\\
3.40082179999445	5.06002491404293\\
4.77363059995696	4.77667373051705\\
7.26964660000522	4.57707006764094\\
9.64086179994047	4.42507464956054\\
11.6376746000024	4.29311535638748\\
13.7592881999444	4.15732466205045\\
16.0057025999995	4.06914434470653\\
18.2521169999382	3.96119206446315\\
20.6233321999898	3.87585558717682\\
23.4937505999114	3.79275988145233\\
25.8649657999631	3.6866923278836\\
28.7353842000011	3.63038038261034\\
31.7306033999193	3.59193022974517\\
34.4762209999608	3.55112138756421\\
37.471440199879	3.50783970120968\\
40.8410617999034	3.48705832576353\\
43.4618785999483	3.47150663045471\\
46.3322969998699	3.45924198930155\\
48.8283129999181	3.44920263184368\\
51.5739305999596	3.44078731625236\\
54.8343514999142	3.43421246261438\\
58.079172299942	3.43070871410766\\
61.4487938999664	3.42738381375213\\
65.9416226999601	3.42115586500824\\
69.6856466998579	3.41520282136783\\
73.3048698998755	3.41137291050434\\
76.5496906999033	3.40750202803593\\
80.0441130999243	3.40377936476731\\
83.1641330999555	3.40011304407417\\
86.7833562998567	3.39617362238874\\
89.9033762997715	3.39203246280603\\
93.5225994997891	3.38768107881226\\
98.1402290997794	3.38303757735749\\
102.38345629978	3.37824903955962\\
105.503476299811	3.37348068468192\\
108.873097899836	3.36890816472437\\
112.61712189985	3.36467269715585\\
116.361145899864	3.36115264062007\\
120.354771499871	3.35828388405347\\
123.599592299783	3.35592761656574\\
127.59321789979	3.35423740657108\\
131.836445099791	3.35285813724321\\
135.830070699798	3.35185901633073\\
139.324493099819	3.35098088372188\\
144.363325399812	3.35034071368936\\
148.856154199806	3.3498944019755\\
153.72338539979	3.34949957309996\\
158.091413399787	3.34911457627257\\
162.833843799774	3.34887004551873\\
167.82587579987	3.34868390103383\\
172.568306199857	3.34849500493736\\
177.685138999834	3.34831739717445\\
182.427569399821	3.34818070539874\\
184.29958139977	3.34807947588771\\
186.670796599821	3.34801538378178\\
189.166812599869	3.34797005648034\\
191.787629399798	3.34789640321076\\
194.283645399846	3.34782319311267\\
196.405258999788	3.34779772290366\\
199.04167589976	3.3477782040913\\
201.662492699805	3.3477373001841\\
203.90890709986	3.34769779913119\\
206.779325499781	3.34766145911694\\
209.27534149983	3.34764663469974\\
211.521755899768	3.34762281998074\\
213.768170299823	3.34759860410849\\
217.808596199844	3.34759098637478\\
220.304612199776	3.34757245283431\\
223.050229799817	3.34755440688097\\
226.045448999852	3.34752638422135\\
229.165468999767	3.34749690343533\\
232.410289799795	3.34747761847885\\
235.155907399836	3.34746738306767\\
238.026325799758	3.34746378661995\\
240.647142599802	3.34746307258206\\
243.657961899764	3.34746145510143\\
247.277185099781	3.34746027348737\\
250.14760349982	3.34745919893054\\
252.518818699871	3.34744635263903\\
255.014834699803	3.34742750661081\\
258.010053899838	3.34742038546519\\
260.256468299776	3.34741633333558\\
262.892885199864	3.34741749672954\\
266.012905199779	3.34741967663356\\
269.257725999807	3.3474164343982\\
272.128144399845	3.34741854199888\\
274.87376199977	3.34742339788904\\
277.868981199805	3.34741899461979\\
280.739399599843	3.34740859388272\\
283.360216399771	3.34740125878164\\
285.731431599823	3.34738625633781\\
289.52225589985	3.34738357686429\\
292.018271899782	3.34737478632206\\
294.51428789983	3.34736818244682\\
297.135104699759	3.34735269056284\\
};
\addplot [color=red,solid,forget plot]
  table[row sep=crcr]{%
0.0725779719550687	1099.7\\
0.124800800113007	428.341302934067\\
0.2496016001096	85.1774040645901\\
0.577203700086102	20.2496617476092\\
1.20120770006906	6.70323948650063\\
2.94841890002135	5.55515768956423\\
4.32122770010028	5.06002491404293\\
6.31804050004575	4.77667373051705\\
8.92325720004737	4.57707006764094\\
11.294472400099	4.42507464956054\\
13.6656876000343	4.29311535638748\\
16.2865044000791	4.15732466205045\\
19.0321220000042	4.06914434470653\\
21.9961410000687	3.96119206446315\\
24.8665594001068	3.87585558717682\\
27.7369778000284	3.79275988145233\\
30.5137956000399	3.6866923278836\\
33.5090148000745	3.63038038261034\\
36.254632400116	3.59193022974517\\
39.000250000041	3.55112138756421\\
41.6210668000858	3.50783970120968\\
44.8658876001136	3.48705832576353\\
48.110708400025	3.47150663045471\\
51.605130800046	3.45924198930155\\
55.3491548000602	3.44920263184368\\
58.7187764000846	3.44078731625236\\
62.3379996001022	3.43421246261438\\
67.252031100099	3.43070871410766\\
71.1208559001097	3.42738381375213\\
75.5356842000037	3.42115586500824\\
78.655704200035	3.41520282136783\\
81.9005250000628	3.41137291050434\\
85.1453458000906	3.40750202803593\\
88.0001641000854	3.40377936476731\\
91.4945865001064	3.40011304407417\\
95.3634113001172	3.39617362238874\\
98.7330329000251	3.39203246280603\\
101.603451300063	3.38768107881226\\
104.848272100091	3.38303757735749\\
107.968292100006	3.37824903955962\\
110.96351130004	3.37348068468192\\
114.426733500091	3.36890816472437\\
118.295558300102	3.36467269715585\\
121.66517990001	3.36115264062007\\
126.579211400007	3.35828388405347\\
130.448036200018	3.35592761656574\\
134.441661800025	3.35423740657108\\
138.560088200029	3.35285813724321\\
142.179311400047	3.35185901633073\\
146.422538600047	3.35098088372188\\
150.416164200054	3.35034071368936\\
154.409789800062	3.3498944019755\\
159.027419400052	3.34949957309996\\
163.395447400049	3.34911457627257\\
167.26427220006	3.34887004551873\\
171.6479003001	3.34868390103383\\
176.515131500084	3.34849500493736\\
181.756765100057	3.34831739717445\\
186.998398700031	3.34818070539874\\
189.494414700079	3.34807947588771\\
191.740829100017	3.34801538378178\\
193.987243500073	3.34797005648034\\
196.608060300117	3.34789640321076\\
199.852881100029	3.34782319311267\\
203.097701900057	3.34779772290366\\
205.609318000032	3.3477782040913\\
207.855732400087	3.3477373001841\\
210.102146800025	3.34769779913119\\
212.800964100054	3.34766145911694\\
215.047378500109	3.34764663469974\\
217.418593700044	3.34762281998074\\
219.914609700092	3.34759860410849\\
222.660227300017	3.34759098637478\\
225.780247300048	3.34757245283431\\
228.75986640004	3.34755440688097\\
231.505484000081	3.34752638422135\\
234.375902400119	3.34749690343533\\
236.996719200048	3.34747761847885\\
240.116739200079	3.34746738306767\\
243.6111616001	3.34746378661995\\
246.980783200008	3.34746307258206\\
249.726400800049	3.34746145510143\\
252.721620000084	3.34746027348737\\
255.592038400006	3.34745919893054\\
258.35325610009	3.34744635263903\\
261.707277600071	3.34742750661081\\
264.577696000109	3.34742038546519\\
267.697716000024	3.34741633333558\\
270.692935200059	3.34741749672954\\
273.563353600097	3.34741967663356\\
276.683373600012	3.3474164343982\\
279.678592800046	3.34741854199888\\
282.673812000081	3.34742339788904\\
285.169828000013	3.34741899461979\\
287.790644800058	3.34740859388272\\
290.286660800106	3.34740125878164\\
293.406680800021	3.34738625633781\\
296.027497600066	3.34738357686429\\
298.273912000004	3.34737478632206\\
301.269131200039	3.34736818244682\\
};
\addplot [color=red,solid,forget plot]
  table[row sep=crcr]{%
0.0725779718873673	1099.7\\
0.124800799996592	428.341302934067\\
0.249601599993184	85.1774040645901\\
0.499203199986368	20.2496617476092\\
1.12320719996933	6.70323948650063\\
2.87041840003803	5.55515768956423\\
4.36802799999714	5.06002491404293\\
5.97483830002602	4.77667373051705\\
8.34605349996127	4.57707006764094\\
10.5924679000163	4.42507464956054\\
13.2132846999448	4.29311535638748\\
15.709300699993	4.15732466205045\\
18.7045199000277	4.06914434470653\\
21.6061385000357	3.96119206446315\\
24.4921570000006	3.87585558717682\\
26.8633722000523	3.79275988145233\\
29.3593881999841	3.6866923278836\\
31.6058026000392	3.63038038261034\\
33.8522169999778	3.59193022974517\\
36.8474362000125	3.55112138756421\\
39.7178546000505	3.50783970120968\\
42.5102724999888	3.48705832576353\\
46.1606959000928	3.47150663045471\\
49.5303175000008	3.45924198930155\\
53.1339406000916	3.44920263184368\\
55.8795582000166	3.44078731625236\\
58.5783755000448	3.43421246261438\\
62.5252008000389	3.43070871410766\\
66.5188264000462	3.42738381375213\\
70.652852899977	3.42115586500824\\
74.3968768999912	3.41520282136783\\
78.0785005000653	3.41137291050434\\
82.4465285000624	3.40750202803593\\
85.6913493000902	3.40377936476731\\
89.4665735000744	3.40011304407417\\
93.2105975000886	3.39617362238874\\
96.8298207001062	3.39203246280603\\
99.7002391001442	3.38768107881226\\
102.945059900172	3.38303757735749\\
106.31468150008	3.37824903955962\\
109.309900700115	3.37348068468192\\
112.429920700146	3.36890816472437\\
115.674741500174	3.36467269715585\\
118.919562300085	3.36115264062007\\
122.663586300099	3.35828388405347\\
126.033207900124	3.35592761656574\\
130.151634300128	3.35423740657108\\
133.895658300142	3.35285813724321\\
138.138885500142	3.35185901633073\\
141.75810870016	3.35098088372188\\
146.812541100197	3.35034071368936\\
150.556565100211	3.3498944019755\\
155.252195200184	3.34949957309996\\
158.996219200199	3.34911457627257\\
163.738649600185	3.34887004551873\\
168.855482400162	3.34868390103383\\
174.221916800132	3.34849500493736\\
178.589944800129	3.34831739717445\\
183.332375200116	3.34818070539874\\
185.828391200164	3.34807947588771\\
187.82520400011	3.34801538378178\\
190.570821600151	3.34797005648034\\
193.815642400179	3.34789640321076\\
196.561260000104	3.34782319311267\\
199.057276000152	3.34779772290366\\
202.208496200154	3.3477782040913\\
205.203715400188	3.3477373001841\\
207.574930600123	3.34769779913119\\
210.070946600172	3.34766145911694\\
212.94136500021	3.34764663469974\\
215.936584200128	3.34762281998074\\
218.931803400163	3.34759860410849\\
221.178217800101	3.34759098637478\\
224.048636200139	3.34757245283431\\
226.419851400191	3.34755440688097\\
229.415070600109	3.34752638422135\\
232.035887400154	3.34749690343533\\
234.906305800192	3.34747761847885\\
238.151126600103	3.34746738306767\\
241.146345800138	3.34746378661995\\
244.141565000173	3.34746307258206\\
247.386385800201	3.34746145510143\\
250.506405800115	3.34746027348737\\
253.12722260016	3.34745919893054\\
255.748039400205	3.34744635263903\\
258.368856200133	3.34742750661081\\
261.488876200165	3.34742038546519\\
264.218893700163	3.34741633333558\\
267.338913700194	3.34741749672954\\
269.834929700126	3.34741967663356\\
272.330945700174	3.3474164343982\\
275.326164900209	3.34741854199888\\
277.96258180018	3.34742339788904\\
280.458597800112	3.34741899461979\\
283.32901620015	3.34740859388272\\
286.324235400185	3.34740125878164\\
289.31945460022	3.34738625633781\\
292.189873000141	3.34738357686429\\
295.185092200176	3.34737478632206\\
298.18031140021	3.34736818244682\\
300.801128200139	3.34735269056284\\
};
\addplot [color=red,solid,forget plot]
  table[row sep=crcr]{%
0	428.341302934067\\
0.2496016001096	85.1774040645901\\
0.374402400106192	20.2496617476092\\
0.87360560009256	6.70323948650063\\
2.74561760004144	5.55515768956423\\
3.99362560000736	5.06002491404293\\
5.49123520008288	4.77667373051705\\
7.86245040001813	4.57707006764094\\
10.483267200063	4.42507464956054\\
12.9792831999948	4.29311535638748\\
15.2256976000499	4.15732466205045\\
18.2209168000845	4.06914434470653\\
21.2161360000027	3.96119206446315\\
24.0865544000408	3.87585558717682\\
26.7073712000856	3.79275988145233\\
29.4529888000106	3.6866923278836\\
32.4482080000453	3.63038038261034\\
35.0690247999737	3.59193022974517\\
38.8130487999879	3.55112138756421\\
42.5570728000021	3.50783970120968\\
45.8018936000299	3.48705832576353\\
48.7971128000645	3.47150663045471\\
51.7923319999827	3.45924198930155\\
56.2383604999632	3.44920263184368\\
59.483181299991	3.44078731625236\\
62.977603700012	3.43421246261438\\
66.7684280000394	3.43070871410766\\
70.4500515999971	3.42738381375213\\
73.6948724000249	3.42115586500824\\
77.5636972000357	3.41520282136783\\
81.6821236000396	3.41137291050434\\
85.5509484000504	3.40750202803593\\
89.1701715999516	3.40377936476731\\
93.0389963999623	3.40011304407417\\
96.2838171999902	3.39617362238874\\
99.9030404000077	3.39203246280603\\
103.522263600025	3.38768107881226\\
107.26628760004	3.38303757735749\\
110.760710000061	3.37824903955962\\
114.130331599968	3.37348068468192\\
117.375152399996	3.36890816472437\\
120.994375600014	3.36467269715585\\
124.363997200038	3.36115264062007\\
127.858419600059	3.35828388405347\\
131.727244399954	3.35592761656574\\
135.596069199964	3.35423740657108\\
139.464893999975	3.35285813724321\\
143.458519599983	3.35185901633073\\
147.82654759998	3.35098088372188\\
151.944973999984	3.35034071368936\\
156.687404399971	3.3498944019755\\
161.429834799957	3.34949957309996\\
165.673061999958	3.34911457627257\\
170.914695600048	3.34887004551873\\
175.157922800048	3.34868390103383\\
179.775552400039	3.34849500493736\\
184.143580400036	3.34831739717445\\
189.260413200012	3.34818070539874\\
191.756429200061	3.34807947588771\\
194.252445199993	3.34801538378178\\
196.748461200041	3.34797005648034\\
199.119676399976	3.34789640321076\\
201.615692400024	3.34782319311267\\
204.111708399956	3.34779772290366\\
206.857325999998	3.3477782040913\\
209.602943600039	3.3477373001841\\
212.348561199964	3.34769779913119\\
215.094178800005	3.34766145911694\\
217.590194800054	3.34764663469974\\
219.961409999989	3.34762281998074\\
222.457426000037	3.34759860410849\\
225.203043599962	3.34759098637478\\
228.073462	3.34757245283431\\
230.694278799929	3.34755440688097\\
233.689497999963	3.34752638422135\\
236.435115600005	3.34749690343533\\
239.679936399916	3.34747761847885\\
242.675155599951	3.34746738306767\\
245.795175599982	3.34746378661995\\
248.66559400002	3.34746307258206\\
251.660813199938	3.34746145510143\\
254.281629999983	3.34746027348737\\
256.902446800028	3.34745919893054\\
259.273661999963	3.34744635263903\\
262.019279600005	3.34742750661081\\
265.139299599919	3.34742038546519\\
268.259319599951	3.34741633333558\\
271.129737999989	3.34741749672954\\
274.124957200023	3.34741967663356\\
276.870574799948	3.3474164343982\\
279.61619239999	3.34741854199888\\
282.736212400021	3.34742339788904\\
285.357029199949	3.34741899461979\\
288.227447599988	3.34740859388272\\
290.723463599919	3.34740125878164\\
293.219479599968	3.34738625633781\\
296.214698800002	3.34738357686429\\
299.20991799992	3.34737478632206\\
301.830734799965	3.34736818244682\\
};
\addplot [color=red,solid,forget plot]
  table[row sep=crcr]{%
0.0725779718873673	1099.7\\
0.124800799996592	428.341302934067\\
0.249601599993184	85.1774040645901\\
0.499203199986368	20.2496617476092\\
1.12320719996933	6.70323948650063\\
2.37121519993525	5.55515768956423\\
3.49442239990458	5.06002491404293\\
4.86723119998351	4.77667373051705\\
7.61284879990853	4.57707006764094\\
9.85926319996361	4.42507464956054\\
12.1056775999023	4.29311535638748\\
14.8512951999437	4.15732466205045\\
17.8465143999783	4.06914434470653\\
20.7169328000164	3.96119206446315\\
23.4625503999414	3.87585558717682\\
26.5825703999726	3.79275988145233\\
29.2033872000175	3.6866923278836\\
32.635409199982	3.63038038261034\\
35.2562259999104	3.59193022974517\\
37.8770427999552	3.55112138756421\\
41.2466643998632	3.50783970120968\\
44.491485199891	3.48705832576353\\
47.7207058998756	3.47150663045471\\
50.8407258999068	3.45924198930155\\
53.7111442999449	3.44920263184368\\
57.2055666998494	3.44078731625236\\
60.6999890998704	3.43421246261438\\
64.1944114998914	3.43070871410766\\
67.813634699909	3.42738381375213\\
71.0740555998636	3.42115586500824\\
74.5684779998846	3.41520282136783\\
78.1877011999022	3.41137291050434\\
81.1829203999368	3.40750202803593\\
84.4277411999647	3.40377936476731\\
87.9221635998692	3.40011304407417\\
91.5413867998868	3.39617362238874\\
95.2698106998578	3.39203246280603\\
98.7642330998788	3.38768107881226\\
102.2586554999	3.38303757735749\\
105.628277099924	3.37824903955962\\
109.122699499945	3.37348068468192\\
112.866723499959	3.36890816472437\\
116.361145899864	3.36467269715585\\
119.730767499888	3.36115264062007\\
123.474791499902	3.35828388405347\\
127.203215399873	3.35592761656574\\
131.072040199884	3.35423740657108\\
135.876870999928	3.35285813724321\\
139.745695799938	3.35185901633073\\
143.817321899929	3.35098088372188\\
148.310150699923	3.35034071368936\\
152.927780299913	3.3498944019755\\
157.420609099907	3.34949957309996\\
161.913437899901	3.34911457627257\\
166.655868299888	3.34887004551873\\
171.148697099881	3.34868390103383\\
175.891127499868	3.34849500493736\\
180.508757099858	3.34831739717445\\
185.375988299842	3.34818070539874\\
187.622402699897	3.34807947588771\\
190.118418699829	3.34801538378178\\
192.364833099884	3.34797005648034\\
195.110450699925	3.34789640321076\\
197.731267499854	3.34782319311267\\
200.383284499869	3.34779772290366\\
203.128902099794	3.3477782040913\\
205.250515699852	3.3477373001841\\
207.496930099791	3.34769779913119\\
209.868145299843	3.34766145911694\\
212.613762899884	3.34764663469974\\
215.608982099802	3.34762281998074\\
218.229798899847	3.34759860410849\\
221.225018099882	3.34759098637478\\
224.469838899793	3.34757245283431\\
228.213862899807	3.34755440688097\\
231.583484499832	3.34752638422135\\
234.578703699866	3.34749690343533\\
237.199520499795	3.34747761847885\\
240.444341299823	3.34746738306767\\
244.687568499823	3.34746378661995\\
248.181990899844	3.34746307258206\\
250.927608499886	3.34746145510143\\
253.922827699804	3.34746027348737\\
257.042847699835	3.34745919893054\\
260.287668499863	3.34744635263903\\
263.407688499894	3.34742750661081\\
266.527708499809	3.34742038546519\\
269.148525299854	3.34741633333558\\
272.143744499888	3.34741749672954\\
274.514959699824	3.34741967663356\\
277.634979699855	3.3474164343982\\
280.380597299896	3.34741854199888\\
282.876613299828	3.34742339788904\\
285.24782849988	3.34741899461979\\
287.743844499812	3.34740859388272\\
290.863864499843	3.34740125878164\\
293.859083699877	3.34738625633781\\
296.854302899796	3.34738357686429\\
299.724721299834	3.34737478632206\\
302.470338899875	3.34736818244682\\
};
\end{axis}
\end{tikzpicture}%

%% file: cmm_solvers.tikz
%
%
\begin{tikzpicture}
\tikzstyle{every node}=[font=\footnotesize]

\begin{axis}[%
width=\figurewidth,
height=\figureheight,
at={(0\figurewidth,0\figureheight)},
scale only axis,
separate axis lines,
every outer x axis line/.append style={black},
every x tick label/.append style={font=\color{black}},
every x tick label/.append style={font=\scriptsize},
xmode=log,
xmin=0.1,
xmax=300,
xminorticks=true,
xlabel={Time (sec)},
every outer y axis line/.append style={black},
every y tick label/.append style={font=\color{black}},
every y tick label/.append style={font=\scriptsize},
ymode=log,
ymin=3,
ymax=2000,
yminorticks=true,
y label style={at={(axis description cs:0.05,0.45)},anchor=north},
ylabel={Cost},
legend style={legend cell align=left,align=left,draw=black}
]
\addplot [color=red,solid]
  table[row sep=crcr]{%
0	2032.73741689515\\
0.124800799996592	1073.28957020101\\
0.795605099992827	458.991076829879\\
0.920405899989419	194.681656944674\\
1.04520669998601	90.6850807775417\\
1.1700074999826	40.3061304900791\\
1.41960909997579	20.0890178086842\\
1.66921069996897	11.1808054598945\\
2.04361309995875	7.19069850605226\\
3.04201949993148	5.66782820192258\\
3.66602350003086	5.11933399766954\\
4.78923070000019	4.8252345375313\\
6.1620394999627	4.64852702545092\\
7.65964909992181	4.51426793347809\\
8.78285630000755	4.37496240451404\\
10.1556650999701	4.23174997420057\\
11.5284738999326	4.12274160338278\\
13.0260835000081	4.04353427351866\\
14.3988922999706	3.97427168881271\\
15.8965018999297	3.9035958963558\\
17.5189123000018	3.84070757648536\\
19.1413226999575	3.79238429718633\\
20.638932300033	3.74674530881433\\
21.886940299999	3.71043968901842\\
23.2597490999615	3.69745240427198\\
24.632557899924	3.66018035620856\\
27.2845748999389	3.64737723249606\\
29.9053916999837	3.63606109059912\\
31.9022044999292	3.6248494519764\\
33.7742164999945	3.61176584138319\\
35.6462284999434	3.58281289644268\\
37.1438381000189	3.54089188853105\\
38.641447699978	3.50075065880141\\
40.2638580999337	3.4757944306317\\
42.0110693000024	3.47649645800782\\
43.5398790999316	3.44785333064897\\
45.2870903000003	3.43602765408147\\
46.7846998999594	3.43130165081405\\
48.1575086999219	3.42773159564281\\
49.4055166998878	3.42479519710969\\
50.7783254999667	3.42187532163882\\
52.260334999999	3.41935525683699\\
53.9919461000245	3.41700835012605\\
55.364754899987	3.41489058981541\\
56.8623644999461	3.4128951321412\\
58.2351733000251	3.41079738221933\\
59.7327828999842	3.40894777997381\\
61.604794899933	3.4070788759977\\
63.6172078000382	3.40524146743271\\
65.489219799987	3.40346653273591\\
68.3908383999951	3.40168227809997\\
70.6372527999338	3.3997147801774\\
72.6652657999657	3.39777370683807\\
74.6620786000276	3.39559587111819\\
76.5340905999765	3.3933345916078\\
78.4061025999254	3.39083528123247\\
80.2781145999907	3.38806449002769\\
82.1501265999395	3.38501173488385\\
84.0221386000048	3.38190468176934\\
85.8941505999537	3.37852290286446\\
87.3917602000292	3.37504121429446\\
89.2637721999781	3.37156263342471\\
90.9485829999903	3.36860889447862\\
93.350998400012	3.36575796137735\\
95.2230103999609	3.36330848186802\\
97.3446240000194	3.36112855808159\\
99.0918351999717	3.35928324031252\\
100.839046399924	3.35790775181932\\
102.336655999999	3.35677695916496\\
104.208667999948	3.35577939217475\\
105.706277600024	3.35482957039407\\
107.32868799998	3.35406289434014\\
108.951098399935	3.35330367394576\\
110.573508800007	3.35278909982069\\
112.648322099936	3.35226046955748\\
114.395533300005	3.35174153532234\\
116.267545299954	3.35136946060769\\
118.014756500022	3.35098013599365\\
119.886768499971	3.35063230645681\\
121.63397970004	3.35037720587971\\
124.005194899975	3.35007268990599\\
126.376410100027	3.34980617501697\\
128.622824499966	3.34957610448911\\
131.118840500014	3.34942689465626\\
133.240454099956	3.34919718530993\\
135.237266900018	3.34901975767605\\
137.73328289995	3.34883195986995\\
140.229298899998	3.34874786181026\\
142.72531489993	3.34859472620133\\
145.221330899978	3.34846131801159\\
147.21814370004	3.3483214246095\\
149.589358899975	3.34820545022409\\
152.085374900023	3.34811033675781\\
154.706191699952	3.34803686831486\\
157.451809299993	3.3479382311224\\
160.072626100038	3.34785149034498\\
162.319040499977	3.34778529534856\\
164.440654100035	3.34772571851134\\
166.562267699977	3.34767287639544\\
169.432686100015	3.34765849071226\\
172.303104500053	3.34762081806492\\
};
\addlegendentry{Trust regions};

\addplot [color=blue,solid]
  table[row sep=crcr]{%
0	2032.73741689515\\
0.405602599959821	471.970758531747\\
0.405602599959821	161.680577090692\\
0.561603599926457	101.267513889198\\
0.639604099909775	52.3918409046649\\
0.717604599893093	37.3123863332648\\
0.795605099876411	33.4683853836707\\
0.873605599859729	22.002602566146\\
0.951606099959463	19.4357278986441\\
1.02960659994278	17.4704379076906\\
1.1076070999261	16.6103540338796\\
1.18560759990942	13.7747141162841\\
1.26360809989274	12.9155508487524\\
1.41960909985937	11.8881930739689\\
1.49760959995911	11.3893723715121\\
1.57561009994242	11.1418129158883\\
1.65361059992574	10.0934046712191\\
1.73161109990906	9.62741360407886\\
1.80961159989238	9.36377402844717\\
1.96561259985901	9.14395768969048\\
2.04361309995875	8.94812644812125\\
2.12161359994207	8.74927239797237\\
2.19961409992538	8.61234972539947\\
2.2776145999087	8.14367866343964\\
2.43361559987534	7.99015251315551\\
2.51161609985866	7.70846297640421\\
2.58961659995839	7.65491621217886\\
2.66761709994171	7.35904257315268\\
2.74561759992503	7.29044109363784\\
2.90161859989166	7.23029635627851\\
2.97961909987498	7.13983105948154\\
3.0576195998583	7.09406514969588\\
3.13562009995803	6.81688731483579\\
3.21362059994135	6.77353035942299\\
3.29162109992467	6.54317557307997\\
3.43202199996449	6.46903189044313\\
3.51002249994781	6.41614529397946\\
3.58802299993113	6.37966034927282\\
3.66602349991445	6.34085577348692\\
3.82202449988108	6.31443003539231\\
3.94682529987767	6.23074890083601\\
4.19642689987086	6.21361066967719\\
4.38362809992395	5.85847631120259\\
4.50842889992055	5.81114147688111\\
4.72683029994369	5.78651206961019\\
4.85163109994028	5.76704132731488\\
4.97643189993687	5.74952829981427\\
5.27283379994333	5.73267989462177\\
5.39763459993992	5.72008217571573\\
5.52243539993651	5.66236300832391\\
5.6472361999331	5.64120469789746\\
5.7720369999297	5.57296688995244\\
5.89683779992629	5.55974219262682\\
6.02163859992288	5.49795593934737\\
6.27124019991606	5.47863671338361\\
6.39604099991266	5.4633388750652\\
6.52084179990925	5.44772093948442\\
6.64564259990584	5.43371408234976\\
6.77044339990243	5.41707781645072\\
6.89524419989903	5.40647821510141\\
7.08244539995212	5.35758005353526\\
7.20724619994871	5.33297908982562\\
7.41004749992862	5.31996854017227\\
7.53484829992522	5.30474060139569\\
7.7844498999184	5.29602429348293\\
7.7844498999184	5.23294302975769\\
8.03405149991158	5.22442924543249\\
8.15885229990818	5.2108822571439\\
8.28365309990477	5.202001467933\\
8.40845389990136	5.19416426816126\\
8.65805549989454	5.18599186035169\\
8.78285629989114	5.17843866760762\\
8.90765709988773	5.1742929879343\\
9.03245789988432	5.0204805167236\\
9.20405899989419	5.01193106911642\\
9.34445989993401	4.9965998264804\\
9.4692606999306	4.9928951628269\\
9.71886229992379	4.93712472834697\\
9.84366309992038	4.92309935932911\\
9.9840639999602	4.91238271898834\\
10.1244648998836	4.90841013549928\\
10.2492656998802	4.87329061655265\\
10.3740664998768	4.86448499161596\\
10.4988672998734	4.85390740214796\\
10.62366809987	4.84691554011198\\
10.8108692999231	4.84396959359234\\
10.9356700999197	4.79099198853572\\
10.9356700999197	4.78202036249486\\
11.0604708999163	4.77849744545189\\
11.1852716999128	4.75691431002982\\
11.3100724999094	4.75004212968676\\
11.5128737998893	4.74687259655887\\
11.6376745998859	4.7373364335535\\
11.9028762999224	4.73406400778568\\
12.0432771999622	4.71730444283489\\
12.1680779999588	4.71090516158997\\
12.2928787999554	4.70842337744044\\
12.5424803999485	4.67399329889492\\
12.6672811999451	4.66695277538447\\
12.7920819999417	4.66417940166055\\
12.9168827999383	4.6607715883141\\
};
\addlegendentry{Conjugate gradients};

\addplot [color=red,solid]
  table[row sep=crcr]{%
0	2032.73741689515\\
0.124800799996592	428.341302934067\\
0.374402399989776	85.1774040645901\\
0.624003999982961	20.2496617476092\\
1.49760959995911	6.70323948650063\\
3.1200199999148	5.55515768956423\\
4.86723119998351	5.06002491404293\\
6.6144423999358	4.77667373051705\\
8.86085679999087	4.57707006764094\\
12.3708792999387	4.42507464956054\\
15.1164968999801	4.29311535638748\\
17.7373136999086	4.15732466205045\\
21.7777395999292	4.06914434470653\\
24.7729587999638	3.96119206446315\\
27.3937756000087	3.87585558717682\\
31.6682029999793	3.79275988145233\\
34.4138205999043	3.6866923278836\\
37.9394431998953	3.63038038261034\\
41.0594631999265	3.59193022974517\\
44.0546823999612	3.55112138756421\\
46.9251007999992	3.50783970120968\\
50.2947223999072	3.48705832576353\\
53.2899415999418	3.47150663045471\\
56.5659625999397	3.45924198930155\\
60.0603849999607	3.44920263184368\\
64.5844139999244	3.44078731625236\\
68.5780395999318	3.43421246261438\\
71.9476611999562	3.43070871410766\\
75.4420835999772	3.42738381375213\\
79.1861075999914	3.42115586500824\\
82.6961300999392	3.41520282136783\\
86.2997531999135	3.41137291050434\\
91.1981845999835	3.40750202803593\\
94.8174078000011	3.40377936476731\\
98.5614317998989	3.40011304407417\\
102.305455799913	3.39617362238874\\
106.049479799927	3.39203246280603\\
110.027505300008	3.38768107881226\\
113.521927699912	3.38303757735749\\
117.14115089993	3.37824903955962\\
120.885174899944	3.37348068468192\\
124.254796499969	3.36890816472437\\
127.874019699986	3.36467269715585\\
132.77245109994	3.36115264062007\\
136.64127589995	3.35828388405347\\
139.761295899982	3.35592761656574\\
143.754921499989	3.35423740657108\\
147.748547099996	3.35285813724321\\
151.991774299997	3.35185901633073\\
156.391002499964	3.35098088372188\\
160.384628099971	3.35034071368936\\
164.752656099969	3.3498944019755\\
169.495086499956	3.34949957309996\\
173.863114499953	3.34911457627257\\
178.979947299929	3.34887004551873\\
184.081179999979	3.34868390103383\\
189.198012799956	3.34849500493736\\
193.69084159995	3.34831739717445\\
198.433271999937	3.34818070539874\\
200.929287999985	3.34807947588771\\
203.425304000033	3.34801538378178\\
205.921319999965	3.34797005648034\\
208.54213680001	3.34789640321076\\
210.663750399952	3.34782319311267\\
213.409367999993	3.34779772290366\\
216.154985600035	3.3477782040913\\
219.150204799953	3.3477373001841\\
221.895822399994	3.34769779913119\\
224.26703759993	3.34766145911694\\
226.388651199988	3.34764663469974\\
230.039074599976	3.34762281998074\\
233.03429380001	3.34759860410849\\
235.779911399935	3.34759098637478\\
240.303940400016	3.34757245283431\\
243.42396039993	3.34755440688097\\
247.323985400028	3.34752638422135\\
250.319204599946	3.34749690343533\\
252.815220599994	3.34747761847885\\
255.311236600042	3.34746738306767\\
258.96166000003	3.34746378661995\\
261.457675999962	3.34746307258206\\
264.312494299957	3.34746145510143\\
267.432514299988	3.34746027348737\\
270.302932700026	3.34745919893054\\
273.547753499937	3.34744635263903\\
276.667773499968	3.34742750661081\\
279.538191900007	3.34742038546519\\
282.814212900004	3.34741633333558\\
286.433436100022	3.34741749672954\\
289.42865529994	3.34741967663356\\
292.423874499975	3.3474164343982\\
295.419093700009	3.34741854199888\\
298.539113700041	3.34742339788904\\
301.659133699955	3.34741899461979\\
};

\addplot [color=red,solid]
  table[row sep=crcr]{%
0	2032.73741689515\\
0.249601599993184	73.8892753652163\\
1.91881229996216	5.87159306251313\\
6.28684029995929	5.12357887438369\\
12.7764819000149	4.84786330445506\\
17.768513899995	4.65600565895172\\
22.0117410999956	4.49252123689054\\
27.2689748000121	4.36482641801224\\
32.8850107999751	4.23068780330047\\
37.7522419999586	4.10865073437239\\
42.2450707999524	4.04827083416549\\
48.1419085999951	3.98939236587832\\
53.7891448000446	3.93295598587854\\
59.5299816000042	3.88124127544347\\
65.0212167999707	3.81116380531632\\
70.6372527999338	3.75444415562029\\
75.7540856000269	3.71912694815264\\
80.4965160000138	3.67923112155247\\
83.3669343999354	3.62971251358888\\
86.36215359997	3.58894638703586\\
89.4821736000013	3.55355568369266\\
92.2277912000427	3.52617786436766\\
95.2230103999609	3.52643002749498\\
98.8422335999785	3.49631355224286\\
102.211855200003	3.48553514447434\\
105.331875200034	3.4796407517013\\
108.202293599956	3.47507323477351\\
112.211519300006	3.47116234436211\\
116.84474900004	3.46783628389403\\
120.339171399944	3.46499230274654\\
123.973994700005	3.46254181365682\\
128.092421100009	3.46035466240032\\
131.96124590002	3.4585449370715\\
135.705269900034	3.45678799012007\\
139.402493600035	3.45508203759968\\
143.006116700009	3.45361121489413\\
146.750140700024	3.45216361879203\\
150.852966999984	3.45074913487242\\
154.846592599992	3.44930359619747\\
158.465815800009	3.44782182965276\\
162.33464060002	3.44629644489215\\
165.704262200044	3.4442386923475\\
168.824282199959	3.44175118834886\\
172.31870459998	3.43889481663481\\
175.937927799998	3.43550942211738\\
179.182748600026	3.43158035221397\\
182.801971800043	3.42666774696513\\
186.374394700048	3.420530616653\\
189.338413699996	3.41332887954141\\
191.959230500041	3.40489566745534\\
194.580047299969	3.39554856248926\\
197.7000673	3.38567719049397\\
201.069688900025	3.3762929207588\\
203.940107299946	3.36853844790031\\
207.060127299977	3.36284961508879\\
210.679350499995	3.35876008565828\\
213.549768900033	3.35569113461777\\
217.293792900047	3.35369642910409\\
222.441825899994	3.35241258376408\\
227.153056100011	3.3514853701005\\
230.897080099909	3.35072647524842\\
235.639510500012	3.35012651958096\\
240.257140100002	3.34965877044154\\
244.375566500006	3.34929361725906\\
249.24279769999	3.34901672212897\\
254.484431300079	3.34875642813841\\
259.351662500063	3.34855462019759\\
261.223674500012	3.34844493855197\\
263.220487300074	3.34832271729046\\
265.591702500009	3.34825893605583\\
267.713316100067	3.34817064183727\\
270.084531300003	3.34810638558986\\
272.705348100048	3.34805493645335\\
274.82696169999	3.347994573304\\
276.948575300048	3.34793186273367\\
279.194989699987	3.34786367961782\\
281.815806500032	3.3478101980093\\
284.561424100073	3.34778940063376\\
287.166640800075	3.34774440570713\\
289.9122584	3.3476890745326\\
292.657876000041	3.34763977971678\\
294.90429039998	3.34759258149292\\
297.774708800018	3.3475763928435\\
300.395525600063	3.34753741974391\\
};

\addplot [color=blue,solid]
  table[row sep=crcr]{%
0	2032.73741689515\\
0.124800799996592	268.083706742637\\
0.249601599993184	77.042262134316\\
0.374402399989776	51.221689045399\\
0.514803300029598	41.2020812470804\\
0.63960410002619	24.5400252205573\\
0.764404900022782	16.6538804954218\\
0.889205700019374	15.0099040898011\\
1.01400650001597	13.4270476590969\\
1.13880730001256	12.4456681401747\\
1.26360810000915	11.1348675768241\\
1.38840890000574	9.85461691023376\\
1.51320970000234	9.5477333820118\\
1.63801049999893	9.23764823986685\\
1.63801049999893	8.96969763563779\\
1.87201199994888	7.96094964090551\\
1.99681279994547	7.36971344625089\\
2.12161359994207	7.20867149493655\\
2.24641439993866	7.16232303959955\\
2.44921569991857	6.50171778898373\\
2.57401649991516	6.43539895854391\\
2.57401649991516	6.39724722771688\\
2.74561759992503	6.24383819953691\\
2.91721869993489	6.18480279695122\\
3.12002000003122	5.92693952123377\\
3.24482080002781	5.83231225939441\\
3.3696216000244	5.75476887926912\\
3.49442240002099	5.70633522132712\\
3.49442240002099	5.66159112532978\\
3.74402400001418	5.63868615057642\\
3.86882480001077	5.33766445919485\\
3.99362560000736	5.31230743274378\\
4.11842640000395	5.25746648989761\\
4.36802799999714	5.24103280438566\\
4.49282879999373	5.21448681738763\\
4.61762959999032	5.20576761786428\\
4.74243039998692	5.06287354327289\\
4.9920319999801	5.04709460730133\\
5.11683279997669	5.03894631258783\\
5.33523419999983	4.93890671339654\\
5.46003499999642	4.92431324013871\\
5.70963659998961	4.89628585956202\\
5.8344373999862	4.88111426108573\\
5.95923819998279	4.8714105682677\\
6.08403899997938	4.86136611354172\\
6.30244040000252	4.85156160415501\\
6.64564260002226	4.8451411880838\\
6.80164359998889	4.83363538890051\\
7.02004500001203	4.82661998702053\\
7.14484580000862	4.75079723231199\\
7.37884729995858	4.74316987986498\\
7.62844889995176	4.73714546346079\\
7.75324969994836	4.73323213159371\\
7.87805049994495	4.72953168200852\\
8.00285129994154	4.70805466020655\\
8.12765209993813	4.66200871543784\\
8.39285379997455	4.65612692877453\\
8.61125519999769	4.63586794768076\\
8.73605599999428	4.62529948862116\\
8.98565759998746	4.61805227347628\\
9.11045839998405	4.61284459661197\\
9.34445989993401	4.60626114375554\\
9.4692606999306	4.56643730077161\\
9.59406149992719	4.55732776830529\\
9.71886229992379	4.55351441699456\\
9.84366309992038	4.54922273816857\\
9.96846389991697	4.48972947373787\\
10.09326470003	4.48398233125793\\
10.3428663000232	4.46954028425544\\
10.4676671000198	4.46321107710747\\
10.6548682999564	4.45568585288908\\
10.779669099953	4.45068501062001\\
11.0448707999894	4.44235531805547\\
11.169671599986	4.43762534945801\\
11.2944723999826	4.43045394200418\\
11.4192731999792	4.42587361470101\\
11.7780754999258	4.40149384458813\\
11.9028762999224	4.3927457076589\\
12.0276770999189	4.38610438373892\\
12.1524778999155	4.38221357667828\\
12.4020795000251	4.37475830942048\\
12.5268803000217	4.36909857570863\\
12.6516811000183	4.3646600724395\\
12.870082499925	4.35857696134909\\
13.0104833999649	4.35425483741802\\
13.2600849999581	4.2983039257033\\
13.6344873999478	4.29010753387323\\
13.7748882999877	4.25429204084323\\
14.0244898999808	4.24665598598701\\
14.1492906999774	4.23238663684741\\
14.274091499974	4.22317452404077\\
14.3988922999706	4.21693686194566\\
14.5236930999672	4.21078939953919\\
14.6484938999638	4.20352148547355\\
14.8668952999869	4.19889903886681\\
14.9916960999835	4.15940945492618\\
15.1164968999801	4.15464512701516\\
15.3660984999733	4.12973087232838\\
15.4908992999699	4.12428983111371\\
15.6157000999665	4.12070902245534\\
15.7405008999631	4.11437743734207\\
};

\addplot [color=blue,solid]
  table[row sep=crcr]{%
0	2032.73741689515\\
0.187201200053096	120.646409703448\\
0.312002000049688	34.1868275478655\\
0.452402899973094	20.5208900761702\\
0.63960410002619	15.2670964831103\\
0.764404900022782	11.394842803642\\
1.01400650001597	9.42284530091366\\
1.13880730001256	8.49390890896326\\
1.26360810000915	8.19626653840538\\
1.38840890000574	7.99352491443113\\
1.71601099998225	7.60915604177555\\
1.96561259997543	6.51667470808538\\
2.10601350001525	6.27723021044955\\
2.34001499996521	6.11213399904356\\
2.58961659995839	6.01001629595605\\
2.7144174000714	5.7788118881694\\
2.94841890002135	5.69582943341553\\
3.83762460004073	5.64580910323506\\
4.10282629996073	5.55406485426135\\
4.22762709995732	5.41995982114828\\
4.47722870006692	5.30327807353923\\
4.85163110005669	5.23220078438872\\
5.10123270004988	5.1807679435867\\
5.44443490006961	5.14117980368228\\
5.69403650006279	5.08400790016293\\
5.94363810005598	5.04360617214974\\
6.08403899997938	4.96496418033506\\
6.20883979997598	4.93172301479726\\
6.45844139996916	4.83944600408112\\
6.58324219996575	4.81599780951446\\
6.8640440000454	4.77509209395468\\
6.98884480004199	4.74867819198768\\
7.23844640003517	4.73086399825705\\
7.48804800002836	4.69897372274373\\
7.76884979999159	4.66655329678003\\
8.01845139998477	4.65056188690043\\
8.26805299997795	4.62789415649062\\
8.39285379997455	4.58483043435299\\
8.79845640005078	4.57048135437736\\
8.92325720004737	4.53250236153457\\
9.26645940006711	4.51913255440056\\
9.62526170001365	4.45448584814108\\
9.87486330000684	4.43587619747987\\
9.99966410000343	4.4198268616609\\
10.3428663000232	4.40342336865147\\
10.4676671000198	4.37573962461557\\
10.7172687000129	4.36273042504349\\
10.8420695000095	4.3458798634958\\
11.0916711000027	4.33162440025339\\
11.3412726999959	4.31838197266906\\
11.4660734999925	4.27895772683695\\
11.7156750999857	4.26480780311556\\
12.0432771999622	4.2476906517742\\
12.1680779999588	4.21026279106531\\
12.2928788000718	4.19711224399726\\
12.4176796000684	4.1490746430428\\
12.542480400065	4.13663569046389\\
12.9012827000115	4.11381980825495\\
13.2600850000745	4.10204214254212\\
14.1492906999774	4.08911255018835\\
14.4768928000703	4.07801659218401\\
14.9604959000135	4.06919310724597\\
15.2100975000067	4.0527562945022\\
15.4752992000431	4.04343159170097\\
15.7249008000363	4.03519230829549\\
16.0213027000427	4.02122914014477\\
16.3957051000325	4.01412999973346\\
16.6453067000257	3.99717969325842\\
16.7701075000223	3.99088635105562\\
17.0197091000155	3.98517895010483\\
17.1445099000121	3.96431066372328\\
17.3941115000052	3.95873852376181\\
17.6437130999984	3.95232441482666\\
17.8933146999916	3.94633310132924\\
18.0181154999882	3.9354911721735\\
18.4237181000644	3.92887649730223\\
18.548518900061	3.9135245811348\\
18.7981205000542	3.90712925724331\\
18.9229213000508	3.9024426094341\\
19.172522900044	3.89574625958729\\
19.4221245000372	3.89008966396246\\
19.5469253000338	3.88483410243979\\
19.6717261000304	3.87444404829472\\
19.9681280000368	3.8698789868658\\
20.0929288000334	3.85601838969942\\
20.21772960003	3.84976205001405\\
20.4673312000232	3.83237575600079\\
20.6857326000463	3.82719714245564\\
20.8105334000429	3.81921746762463\\
20.9977345999796	3.81041559494104\\
21.2473361999728	3.80534284828468\\
21.4657375999959	3.78643556821546\\
21.5905383999925	3.77701055132947\\
21.8557401000289	3.7678377972146\\
21.9805409000255	3.7624539579643\\
22.2301425000187	3.75659659045353\\
22.3705434000585	3.74697323199599\\
22.7137456000783	3.74456681996628\\
22.8385464000748	3.74239619097632\\
23.088148000068	3.73254211174993\\
23.2129488000646	3.73127626657561\\
};

\node[fill=white, right, align=left, inner sep=0mm, text=red]
at (axis cs:0.212416838182301,763.907784504421,17.3205080756888) {\footnotesize 2};
\node[fill=white, right, align=left, inner sep=0mm, text=red]
at (axis cs:0.226380340952145,178.427581267205,17.3205080756888) {\footnotesize 3};
\node[fill=white, right, align=left, inner sep=0mm, text=red]
at (axis cs:0.368824936656089,39.4903722957669,17.3205080756888) {\footnotesize 5};
\node[fill=white, right, align=left, inner sep=0mm, text=blue]
at (axis cs:0.465803518792816,122.382493696819,17.3205080756888) {\footnotesize 2};
\node[fill=white, right, align=left, inner sep=0mm, text=blue]
at (axis cs:0.144974067037263,188.301509353101,17.3205080756888) {\footnotesize 3};
\node[fill=white, right, align=left, inner sep=0mm, text=blue]
at (axis cs:0.393070133358267,21.8364123397724,17.3205080756888) {\footnotesize 5};
\end{axis}
\end{tikzpicture}%

%% file: cmm_scaling.tikz
%
%
\begin{tikzpicture}
\tikzstyle{every node}=[font=\footnotesize]

\begin{axis}[%
width=\figurewidth,
height=\figureheight,
at={(0\figurewidth,0\figureheight)},
scale only axis,
separate axis lines,
every outer x axis line/.append style={black},
every x tick label/.append style={font=\color{black}},
every x tick label/.append style={font=\scriptsize},
xmin=500,
xmax=5000,
xlabel={Problem size n},
every outer y axis line/.append style={black},
every y tick label/.append style={font=\color{black}},
every y tick label/.append style={font=\scriptsize},
ymode=log,
ymin=0.1,
ymax=100,
yminorticks=true,
y label style={at={(axis description cs:0.05,0.45)},anchor=north},
ylabel={Time/iter (sec)},
legend style={at={(0.03,0.97)},anchor=north west,legend cell align=left,align=left,draw=black}
]
\addplot [color=blue,solid]
  table[row sep=crcr]{%
500	0.0942246040003374\\
1000	0.427442739999387\\
1500	1.3643847459997\\
2000	2.845770242\\
2500	4.61388557599857\\
3000	7.49959207399981\\
3500	10.6168040559976\\
4000	14.8319510759995\\
4500	19.8361511539994\\
5000	25.9566943879984\\
};
\addlegendentry{Lai \& Osher};

\addplot [color=black!50!green,solid]
  table[row sep=crcr]{%
500	0.100152642000467\\
1000	0.480171078001149\\
1500	1.38528887999943\\
2000	2.81269803000148\\
2500	4.63759772800142\\
3000	7.43407165399985\\
3500	10.5915318939998\\
4000	14.8163509759982\\
4500	19.9260077300016\\
5000	26.353248930003\\
};
\addlegendentry{Neumann};

\addplot [color=red,solid,line width=1.5pt]
  table[row sep=crcr]{%
500	0.144144923998974\\
1000	0.176905133998953\\
1500	0.290785864000209\\
2000	0.283921819999814\\
2500	0.383762459999416\\
3000	0.336026153999846\\
3500	0.360986313999165\\
4000	0.453026904000435\\
4500	0.451154891999904\\
5000	0.604035871997476\\
};
\addlegendentry{MADMM};

\end{axis}
\end{tikzpicture}%

%% file: corresp-tosca-kimscurves.tikz
%
%
\begin{tikzpicture}
\tikzstyle{every node}=[font=\footnotesize]

\begin{axis}[%
width=\figurewidth,
height=\figureheight,
at={(0\figurewidth,0\figureheight)},
scale only axis,
separate axis lines,
every outer x axis line/.append style={black},
every x tick label/.append style={font=\color{black}},
every x tick label/.append style={font=\scriptsize},
xmin=0,
xmax=0.25,
xlabel={\% geodesic diameter},
every outer y axis line/.append style={black},
every y tick label/.append style={font=\color{black}},
every y tick label/.append style={font=\scriptsize},
ymin=0,
ymax=1,
y label style={at={(axis description cs:0.05,0.45)},anchor=north},
ylabel={\% of correspondence},
legend style={at={(0.55,0.3)},anchor=north west,legend cell align=left,align=left,draw=black}
]
\addplot [color=blue,solid]
  table[row sep=crcr]{%
0	0.016875\\
0.005	0.0289583333333333\\
0.01	0.053125\\
0.015	0.0795833333333333\\
0.02	0.105208333333333\\
0.025	0.1325\\
0.03	0.162083333333333\\
0.035	0.193958333333333\\
0.04	0.224166666666667\\
0.045	0.251458333333333\\
0.05	0.276041666666667\\
0.055	0.299791666666667\\
0.06	0.321666666666667\\
0.065	0.345416666666667\\
0.07	0.3675\\
0.075	0.389791666666667\\
0.08	0.409791666666667\\
0.085	0.43\\
0.09	0.44875\\
0.095	0.463125\\
0.1	0.475\\
0.105	0.48625\\
0.11	0.497916666666667\\
0.115	0.507708333333333\\
0.12	0.518333333333333\\
0.125	0.526041666666667\\
0.13	0.533333333333333\\
0.135	0.54125\\
0.14	0.5475\\
0.145	0.553125\\
0.15	0.558333333333333\\
0.155	0.561458333333333\\
0.16	0.566666666666667\\
0.165	0.570208333333333\\
0.17	0.574375\\
0.175	0.578958333333333\\
0.18	0.584791666666667\\
0.185	0.588333333333333\\
0.19	0.592708333333334\\
0.195	0.597083333333333\\
0.2	0.600625\\
0.205	0.603333333333333\\
0.21	0.605833333333333\\
0.215	0.610833333333333\\
0.22	0.614375\\
0.225	0.616458333333333\\
0.23	0.619791666666667\\
0.235	0.621666666666667\\
0.24	0.624375\\
0.245	0.626875\\
0.25	0.629791666666667\\
};
\addlegendentry{ LS };

\addplot [color=red,solid,line width=2.0pt]
  table[row sep=crcr]{%
0	0.0816666666666667\\
0.005	0.117083333333333\\
0.01	0.198125\\
0.015	0.288541666666667\\
0.02	0.373333333333333\\
0.025	0.447083333333333\\
0.03	0.510833333333333\\
0.035	0.565833333333334\\
0.04	0.617083333333333\\
0.045	0.658958333333333\\
0.05	0.703333333333333\\
0.055	0.740625\\
0.06	0.779375\\
0.065	0.808333333333333\\
0.07	0.834583333333333\\
0.075	0.855\\
0.08	0.875833333333333\\
0.085	0.896875\\
0.09	0.91125\\
0.095	0.927291666666667\\
0.1	0.936666666666667\\
0.105	0.945833333333333\\
0.11	0.954791666666667\\
0.115	0.961875\\
0.12	0.970416666666667\\
0.125	0.978541666666667\\
0.13	0.98125\\
0.135	0.983125\\
0.14	0.985\\
0.145	0.987083333333333\\
0.15	0.990208333333333\\
0.155	0.991041666666667\\
0.16	0.992083333333333\\
0.165	0.992708333333333\\
0.17	0.994375\\
0.175	0.99625\\
0.18	0.997291666666667\\
0.185	0.998125\\
0.19	0.998541666666667\\
0.195	0.999166666666667\\
0.2	0.999166666666667\\
0.205	0.999375\\
0.21	0.999583333333333\\
0.215	0.999583333333333\\
0.22	0.999791666666667\\
0.225	1\\
0.23	1\\
0.235	1\\
0.24	1\\
0.245	1\\
0.25	1\\
};
\addlegendentry{ MADMM};

\end{axis}
\end{tikzpicture}%

%% file: mds_scale.tikz
%
%
\begin{tikzpicture}
\tikzstyle{every node}=[font=\footnotesize]

\begin{axis}[%
width=\figurewidth,
height=\figureheight,
at={(0\figurewidth,0\figureheight)},
scale only axis,
separate axis lines,
every outer x axis line/.append style={black},
every tick label/.append style={font=\scriptsize},
xmin=0,
xmax=1000,
xlabel={Problem size n},
every outer y axis line/.append style={black},
every y tick label/.append style={font=\color{black}},
every y tick label/.append style={font=\scriptsize},
ymode=log,
ymin=0.001,
ymax=1000,
yminorticks=true,
y label style={at={(axis description cs:0.025,0.45)},anchor=north},
ylabel={Time/iter (sec)},
legend style={at={(0.975,0.025)},anchor=south east,legend cell align=left,align=left,draw=black}
]
\addplot [color=black!50!green,solid]
  table[row sep=crcr]{%
40	0.794529231034489\\
50	5.49279521000002\\
60	12.4717924468749\\
70	29.4770575257143\\
80	61.056261654054\\
90	121.54817915\\
};
\addlegendentry{SDP};

\addplot [color=blue,solid]
  table[row sep=crcr]{%
30	0.00756368484827831\\
40	0.00882429899024598\\
50	0.0266304737370875\\
60	0.0255274363638212\\
70	0.0434911878787762\\
80	0.0354547727273335\\
90	0.0592488646460227\\
100	0.0657095121210347\\
150	0.106737526314407\\
200	0.256991121054039\\
250	0.407335944441406\\
300	0.587603766665173\\
400	1.32600850000098\\
500	5.30750068888948\\
600	8.10685196666357\\
700	10.1781985777779\\
800	12.7712818666672\\
900	18.8691876222211\\
1000	24.1142879111123\\
};
\addlegendentry{Subgradient};

\addplot [color=red,solid,line width=1.5pt]
  table[row sep=crcr]{%
30	0.181795224752554\\
40	0.169129797029839\\
50	0.193225000990154\\
60	0.241260952474711\\
70	0.243423342575102\\
80	0.297637551484657\\
90	0.262575940594381\\
100	0.283118646534847\\
150	0.309030552380829\\
200	0.434574214284242\\
250	0.60840390000323\\
300	1.3174993545452\\
400	2.10317711817333\\
500	3.06471055454511\\
600	4.44744669090935\\
700	6.23720361817819\\
800	7.71354035453194\\
900	9.72028049091618\\
1000	12.034768054537\\
};
\addlegendentry{MADMM};

\end{axis}
\end{tikzpicture}%

%% file: mds_convergence.tikz
%
%
\begin{tikzpicture}
\tikzstyle{every node}=[font=\footnotesize]

\begin{axis}[%
width=\figurewidth,
height=\figureheight,
at={(0\figurewidth,0\figureheight)},
scale only axis,
clip=false,
separate axis lines,
every outer x axis line/.append style={black},
every x tick label/.append style={font=\scriptsize},
xmin=0,
xmax=100,
xlabel={Time (sec)},
every outer y axis line/.append style={black},
every y tick label/.append style={font=\color{black}},
every y tick label/.append style={font=\scriptsize},
ymode=log,
ymin=1900,
ymax=10621,
yminorticks=true,
y label style={at={(axis description cs:0.05,0.45)},anchor=north},
ylabel={Stress},
legend style={legend cell align=left,align=left,draw=black}
]

\addplot [color=blue,solid]
  table[row sep=crcr]{%
0	2618.96016731413\\
1.46640939998906	2609.99106874207\\
3.13562010001624	2603.75681550606\\
4.88283130002674	2598.73343207833\\
6.24004000000423	2594.43301929368\\
8.7516560999793	2590.62428437997\\
10.5456675999449	2587.17525501944\\
12.3396790999686	2584.00632372257\\
14.5548932999372	2581.05918809645\\
16.3801049999311	2578.29980100447\\
19.0321219998878	2575.7004707371\\
20.9665343998931	2573.23858743047\\
22.7449457999319	2570.89452131263\\
25.272161999892	2568.65580397149\\
27.2065743998974	2566.51294456608\\
29.7181904999306	2564.45511622213\\
31.8554041999159	2562.47467011874\\
35.365426699922	2560.56226699552\\
36.816235999926	2558.71336774483\\
38.6414476999198	2556.92314419494\\
40.326258499932	2555.18778864119\\
42.2606708999374	2553.50175909459\\
45.505491699907	2551.86087088375\\
48.3447098999168	2550.26399329745\\
50.1699215999106	2548.70815420517\\
52.2759350999258	2547.19081807993\\
54.974752399954	2545.70895663164\\
56.8467643999611	2544.2614576108\\
58.2195731999236	2542.84811584026\\
61.2147923999	2541.46763625579\\
63.2116051999037	2540.11741154676\\
64.4128128999146	2538.79465314279\\
66.6748273999547	2537.4976912108\\
68.4844389999635	2536.22627735658\\
70.6060525999637	2534.98000904859\\
72.7120660999208	2533.75758627867\\
74.4904774999595	2532.55838326954\\
76.892892899923	2531.38101655507\\
79.7633112999611	2530.22467915436\\
82.2749273999361	2529.08853697979\\
84.3341405999381	2527.97086290138\\
86.6117551999632	2526.87113581767\\
88.1093647999805	2525.78843678145\\
89.8097756999778	2524.72441915064\\
92.3993922999944	2523.67715461194\\
96.1590163999936	2522.6453661478\\
98.9046339999768	2521.6287052509\\
101.213448799972	2520.62693141292\\
};
\addlegendentry{Subgradient};

\addplot [color=blue,solid,forget plot]
  table[row sep=crcr]{%
0	2618.96016731413\\
2.18401399999857	2533.64994859633\\
5.28843389998656	2485.59743460201\\
7.72204949997831	2452.68521467554\\
10.6860684999847	2427.90480160712\\
12.4800800000085	2408.10233249813\\
15.0228963000118	2391.71218498936\\
17.2069103000104	2377.86799644519\\
18.8605208999943	2365.96815059045\\
19.812127000012	2355.55439240348\\
21.6373387000058	2346.34237034752\\
23.1037480999948	2338.08666798479\\
24.5077570999856	2330.65233597079\\
28.2517810999998	2323.9410587238\\
30.1237931000069	2317.83730616278\\
31.4810018000426	2312.2620670873\\
34.6322220000438	2307.16384600031\\
36.2702325000428	2302.48007185844\\
38.8910493000294	2298.15624492043\\
41.6366669000126	2294.14449290011\\
43.2746774000116	2290.42845016002\\
45.9890948000248	2286.9724832769\\
49.3743165000342	2283.73647252544\\
51.2151283000712	2280.69641090599\\
53.3367419000715	2277.84606089498\\
56.4255617000745	2275.18537140815\\
58.4067744000931	2272.68559934911\\
61.7919961001026	2270.33754168643\\
64.5376137000858	2268.11622257364\\
66.3784255001228	2266.01734325578\\
68.5156392001081	2264.02563934786\\
70.4656517000985	2262.12931121828\\
71.807260300091	2260.32567599894\\
73.24246950011	2258.605673644\\
74.896080100094	2256.96843830569\\
76.1128879000898	2255.40670517736\\
77.6572978001204	2253.91352644347\\
79.3733088001609	2252.48274654072\\
81.7445240001543	2251.1101562531\\
83.9285380001529	2249.79032170358\\
85.2545465001604	2248.52272530869\\
86.9081571001443	2247.30239905563\\
89.8877762001357	2246.12653685994\\
92.9297957001836	2244.99320891455\\
94.1934038001928	2243.90323623326\\
96.8766210001777	2242.85593294475\\
100.027841200179	2241.84843898421\\
};
\addplot [color=blue,solid,forget plot]
  table[row sep=crcr]{%
0	2618.96016731413\\
1.71601099998225	2287.35123048519\\
3.44762210000772	2218.18631137904\\
10.6704683999997	2196.26480572738\\
12.7764819000149	2187.29241782287\\
14.5704933999805	2180.96680531143\\
17.5813126999419	2175.88712398853\\
19.1725228999276	2171.50612397666\\
20.6077320998884	2167.59044530107\\
23.6029512998648	2164.05965841355\\
26.4733696998446	2160.85405725899\\
28.5481829998316	2157.9355919816\\
31.5590022997931	2155.22053529526\\
34.647822099796	2152.69747001891\\
36.7694356997963	2150.31827647808\\
39.6242539997911	2148.05927026793\\
41.636666899838	2145.93096519726\\
44.3822844998213	2143.89498552128\\
47.1435021998477	2141.95797107441\\
48.8127128998167	2140.10682594346\\
52.1355341998278	2138.3414298515\\
54.3663484998397	2136.64200904323\\
56.5659625998233	2135.00720912484\\
59.5767818998429	2133.43740445342\\
61.5423944998183	2131.93170123156\\
63.6016076998203	2130.47183980547\\
65.3020185998757	2129.04661563253\\
66.4564259998733	2127.66073202841\\
69.2800440998399	2126.31252711078\\
71.9944614998531	2125.00311024296\\
73.8976736998884	2123.7410887978\\
76.5652907998883	2122.51566817667\\
79.4981095999246	2121.32099205756\\
81.9473252999014	2120.14598814758\\
85.0049448998761	2118.99292465354\\
86.7833562998567	2117.86338310926\\
88.9673702998552	2116.76541732212\\
90.4961800998426	2115.70012478461\\
93.1013967998442	2114.65280762835\\
96.5334187998669	2113.63148852764\\
98.4522310998873	2112.62281421271\\
100.823446299939	2111.62856242012\\
};
\addplot [color=blue,solid,forget plot]
  table[row sep=crcr]{%
0	2618.96016731413\\
1.49760960007552	10620.6392174649\\
3.91562510002404	7401.9498105102\\
6.0216386000393	6021.83639768997\\
7.50364810007159	5226.27813648095\\
10.3740665001096	4701.60479698563\\
12.1524779001484	4327.55106759418\\
13.4004859001143	4046.51413475417\\
14.8200950000901	3827.24488501899\\
16.5985064001288	3651.11243915917\\
19.1881230002036	3506.47879409923\\
21.7933397002053	3385.68824062322\\
23.5873512001708	3283.2702484714\\
26.0833672002191	3195.32760801892\\
28.6417836002074	3119.0832427944\\
30.6229963002261	3052.31976719718\\
31.8866044002352	2993.38911085011\\
34.1954192002304	2941.03515265207\\
36.1922320001759	2894.31146778907\\
37.455840100185	2852.35060670295\\
39.8426554002799	2814.45207925376\\
41.995469200192	2780.03476130831\\
44.5538856001804	2748.66593329372\\
47.4243040002184	2719.97712891885\\
49.218315500184	2693.69448433337\\
50.7471253001131	2669.5761491505\\
52.6503375001485	2647.36780168905\\
54.3975487001007	2626.86474463309\\
57.1587664001854	2607.84383717507\\
59.5143815001938	2432.59333333588\\
61.5111943001393	2440.53976236031\\
62.7748024001485	2427.39811687232\\
65.0056167001603	2387.57145216609\\
67.0804300002055	2389.14958006668\\
69.9664485001704	2377.59161711574\\
72.8368669002084	2339.2774939643\\
74.5840781001607	2336.70688581095\\
77.3608959001722	2302.97192194851\\
79.170507500181	2296.72755119583\\
81.229720700183	2269.43288223966\\
82.4621286002221	2258.99714247849\\
84.6149424001342	2236.90298672253\\
86.3777537001297	2226.11109157815\\
87.8129629001487	2207.06618638159\\
89.6069744001143	2193.69963468895\\
91.4633863000199	2177.23827272307\\
92.7269944000291	2165.80732545039\\
95.2854108000174	2152.52533058227\\
97.1262226001127	2138.77403016528\\
98.4522311000619	2128.94296191835\\
100.63624510006	2113.61378531867\\
};

\node[fill=white, right, align=left, inner sep=0mm, text=blue]
at (axis cs:50.036866359447,2123.80586769058,17.3205080756888) {$\text{\footnotesize 10}^{\text{\tiny -3}}$};
\node[fill=white, right, align=left, inner sep=0mm, text=blue]
at (axis cs:37.8064516129032,2340.12768684324,17.3205080756888) {$\text{\footnotesize 10}^{\text{\tiny -4}}$};
\node[fill=white, right, align=left, inner sep=0mm, text=blue]
at (axis cs:23.8064516129032,2592.49082361235,17.3205080756888) {$\text{\footnotesize 10}^{\text{\tiny -5}}$};
\node[fill=white, right, align=left, inner sep=0mm, text=blue]
at (axis cs:13.8064516129032,3592.49082361235,17.3205080756888) {$\text{\footnotesize 10}^{\text{\tiny -2}}$};

\addplot [color=red,solid,line width=1.5pt]
  table[row sep=crcr]{%
0	2618.96016731413\\
4.35242790001212	2269.47169224056\\
7.1916461000219	2149.92284769474\\
10.6392682000296	2027.8518403837\\
14.3208917999873	2005.49532390871\\
17.9713151999749	1970.14769570818\\
21.2473361999728	1974.75975456771\\
24.382956299989	1959.98478971924\\
27.5809767999453	1944.39871283676\\
31.2782004999463	1936.03144373094\\
35.1002249999437	1932.92129478192\\
38.9690497999545	1929.06954629116\\
44.1794831999578	1927.15488556595\\
47.9859075999702	1925.46244770615\\
52.5411367999623	1924.71154598252\\
54.7719510999741	1924.71154598252\\
58.1415726999403	1922.86550490321\\
63.5392072999384	1922.12100580324\\
65.7700215999503	1922.12100580324\\
71.3080570999882	1921.14705283044\\
74.0224744999432	1921.14705283044\\
77.6416976999026	1921.02873375354\\
82.5401290999143	1919.67312265764\\
84.7709433999262	1919.67312265764\\
88.4681670999271	1920.98681592127\\
93.0857966999174	1918.92151201724\\
95.2230103999027	1918.92151201724\\
97.4850248999428	1918.92151201724\\
102.243055399915	1919.90420705236\\
};
\addlegendentry{MADMM};

\end{axis}
\end{tikzpicture}%

%% file: conclusions.tex
\section{Discussion and Conclusions}

We presented MADMM, a generic algorithm for optimization of non-smooth functions with manifold constraints, and showed that it can be efficiently used in many important problems from the domains of machine learning, computer vision and pattern recognition, and data analysis. Among the key advantages of our method is its remarkable simplicity and lack of parameters to tune - in all our experiments, it worked entirely out-of-the-box. 
We believe that MADMM will be very useful in many other applications in this community. 
%
In our experiments, we observed that MADMM converged independently of the initialization; a theoretical study of convergence properties is an important future direction.  
%
%